\documentclass[final,5p,times,twocolumn,10pt]{elsarticle}
\usepackage{amssymb}
\usepackage{epsfig}
\usepackage{subfigure}
\usepackage{amsmath}
\usepackage{amsfonts}
\usepackage{multicol}
\usepackage[usenames]{color}
\usepackage{booktabs}
\usepackage{changepage}

 \usepackage{amsthm}

\newtheorem{theorem}{Theorem}[section]

\newtheorem{proposition}[theorem]{Proposition}

\newtheorem{definition}[theorem]{Definition}

\begin{document}

\begin{frontmatter}



\title{Characterizing two-timescale nonlinear dynamics using finite-time {L}yapunov
        exponents and vectors\tnoteref{t1}}


\author[UCI]{K.D. Mease\corref{cor1}\fnref{fn1}}\ead{kmease@uci.edu}
\author[UCI]{U.~Topcu\fnref{fn2}}\ead{ utopcu@seas.upenn.edu}
\author[UCI]{E. Aykutlu{\u{g}}\fnref{fn2}}\ead{eaykutlu@uci.edu}
\author[UCI]{M. Maggia\fnref{fn2}}\ead{mmaggia@uci.edu}

\address[UCI]{Department of Mechanical and Aerospace Engineering,
University of California,~Irvine,~CA~92697~USA}
%\address[caltech]{Control and Dynamical Systems MC 107-81
%California Institute of Technology 1200 E. California Blvd Pasadena~CA~91125~USA}
	\cortext[cor1]{Corresponding author}
	\fntext[fn1]{Professor}
	\fntext[fn2]{Graduate student researcher. U. Topcu is currently a Research Assistant Professor at the University
of Pennsylvania. E. Aykutlu{\u{g}} is currently a Postdoctoral Scholar in Earth System Science at U.C. Irvine.}

\tnotetext[t1]{This material is based upon work supported by the National Science Foundation under Grants CMMI-0010085 and CMMI-1069331.}

\begin{abstract}
Finite-time Lyapunov exponents and vectors are used to define and diagnose boundary-layer type,  two-timescale behavior in the tangent linear dynamics and to determine the associated manifold structure in the flow of a finite-dimensional nonlinear autonomous dynamical system. Two-timescale behavior is characterized by a slow-fast splitting of the tangent bundle for a state space region. The slow-fast splitting is defined using finite-time Lyapunov exponents and
vectors, guided by the asymptotic theory of partially hyperbolic sets, with important modifications for the finite-time case; for example, finite-time Lyapunov analysis relies more heavily on the Lyapunov vectors due to their relatively fast convergence compared to that of the corresponding exponents. The splitting is used to locate points on normally hyperbolic center manifolds. Determining manifolds from tangent bundle structure is more generally applicable than approaches, such as the singular perturbation method, that require special normal forms or other {\em a priori} knowledge. The use, features, and accuracy of the approach are illustrated via several detailed examples.

\end{abstract}

\end{frontmatter}


\section{Introduction}
The flow of a finite-dimensional
autonomous nonlinear dynamical system
with multiple timescales may have manifold structure. Characterizing this structure can facilitate simplified analysis
and computation, and lead to greater
understanding of the system behavior. The relevant timescales are most generally in the linear variational dynamics, i.e., tangent
linear dynamics. Our objective is to diagnose two-timescale behavior in tangent linear dynamics with
slow dynamics and both stable and unstable fast dynamics, and to compute the associated manifold structure in the flow
of the nonlinear system. Because the
intent is to analyze finite-time behavior, we first define two-timescale behavior in this context. Though we only
directly
consider two timescales and normally hyperbolic center manifolds in this paper, the discussion and results
are relevant to systems with more than two timescales and also to additional manifold structure, such as the center-stable and center-unstable manifolds relevant to the solution of certain boundary-value problems \cite{ACC09, Guck09,
Rao99, Kopell}. We do not consider
systems with persistent fast oscillations.

Many of the methods available for computing invariant manifolds (i) operate off the linear structure at an equilibrium
point or a periodic orbit \cite{Krauskopf}, or (ii) require \textit{a priori} knowledge of system coordinates adapted to
the manifold
structure, e.g. \cite{Zagaris,Roussel}, or (iii) require \textit{a priori} knowledge of a manifold that can be
analytically or numerically continued to the manifold of interest, e.g. \cite{Broer,DieciOrthog}.
Our particular application context  is flight guidance and control, and our motivation comes from the notable successes
of the singular perturbation method \cite{petar,Omalley} in providing insight and facilitating solution approximation
with reduced-order models
\cite{Naidu}. Geometric singular perturbation theory \cite{Fenichel,Jones} clarifies the manifold structure in the flow
associated with two-timescale behavior.  The singular perturbation method is one means of obtaining the manifold
structure, but it requires a special coordinate representation, i.e., normal form, with a small parameter, such that the
manifold structure for the parameter value of interest can be obtained via matched asymptotic expansions. The singular
perturbation method can be viewed as an analytical continuation method; however there
is no general systematic method of obtaining the required normal form.

The situation of interest is when two-timescale behavior is suspected in a region of state space, perhaps based on
simulation experience, and one wants a means of diagnosing whether or not there are two (or more) disparate timescales
and, if there are, a means of characterizing the associated flow structure. In addition to requiring methodology that
works away from equilibria and periodic orbits and does not require the singularly perturbed normal form, there is the
challenge that, for our target applications, the methodology must be effective when only finite-time behavior is
considered. The approach addressed in this paper, which we refer to as finite-time Lyapunov analysis (FTLA), uses
finite-time Lyapunov exponents (FTLEs) and the associated vectors
(FTLVs), to diagnose two-timescale behavior and characterize
the associated tangent bundle structure, and then uses invariance-based orthogonality conditions to locate and compute the
associated manifold structure. Orthogonality conditions are used in the intrinsic low-dimensional manifold (ILDM) method
\cite{MaasPope} in the chemical kinetics context to compute a slow manifold, but the tangent bundle structure is
determined by a means other than FTLA. Orthogonality conditions are also used for the computation of invariant manifolds
in \cite{DieciOrthog}, but the tangent bundle structure is derived from a known neighboring manifold in a numerical
continuation scheme.

FTLA is used in different ways in several application contexts. The body of work (e.g.,
\cite{Wiggins,Haller11,Sandstede,Marsden}) on characterizing finite-time manifold structure in time-dependent velocity
fields has connections with our work, though the target is co-dimension one manifolds that separate the flow and not
two-timescale behavior. In particular the maximum FTLE field is used to determine Lagrangian coherent structures in
fluid flows with time-dependent velocity fields, e.g., \cite{Doerner,Haller,Marsden,Moser} and to assess the
stability of orbits in celestial mechanics \cite{Froeschle,Villac}. FTLA is used to identify the fastest growing
direction(s) of initialization errors in weather predictability theory \cite{Palmer,Legras,Samelson,Yoden}.

FTLA is applied to systems with slow-fast behavior in \cite{adrover06,JGCD03,AMCpaper,Pope}. In
\cite{JGCD03}, Lyapunov analysis is proposed as a means of diagnosing timescales and suggesting adapted coordinates as an alternative to the singular perturbation approach. In \cite{AMCpaper}, the
ILDM and computational singular perturbation (CSP) \cite{Lam1,Lam2} methods for slow-fast behavior are interpreted
geometrically using Fenichel theory and the idea of using FTLE/Vs to improve the ILDM method is proposed. In
\cite{adrover06} Lyapunov analysis is applied to periodic and chaotic attractors, as well as slow manifolds, and an approach for computing FTLVs is developed. Lorenz
\cite{Lorenz} seems to have been the first to use FTLA to analyze a chaotic attractor. In \cite{Pope}, FTLA is used to identify the
dimension of the attracting slow manifold along a trajectory. The application of FTLA to the solution of two-timescale
boundary value problems related to optimal control is discussed in \cite{ACC09}.

The main contribution of the present paper is to extend FTLA to the diagnosis and computation of normally hyperbolic center manifolds. Because the finite time is limited, it is crucial to define the tangent bundle splitting of interest in the fastest converging way and to clarify the finite-time required to accurately approximate the invariant tangent bundle splitting. Guided by the theory of partially hyperbolic sets \cite{PesinHB}, a finite-time two-timescale set is defined, requiring spatial and temporal uniformity of the spectral gap between the
slow and fast FTLEs. A \textit{fast stable--slow--fast unstable} tangent bundle splitting is specified in terms of the
FTLVs. The size of the spectral gap dictates the rate of exponential
convergence of the tangent bundle splitting toward the desired invariant splitting, providing a guideline for how large the finite-time needs to be. We account for both fast stable and fast unstable behavior
and provide orthogonality conditions for approximately computing points on normally hyperbolic center manifolds, whereas previous
work on FTLA, with the exception of \cite{ACC09}, considered only attracting or repelling center manifolds. Several detailed examples are presented to illustrate and
clarify the approach, and to demonstrate its feasibility and effectiveness in locating and approximating invariant center manifolds.

  The paper is organized as follows.  In Section $\ref{setup}$, we specify the
dynamical system to be considered and recall some definitions from
geometry. Section $\ref{overview}$ provides an overview of the approach and supplements the introduction with background
and perspective required to understand the goals and contributions of the present work as well as relations to other
work. Section $\ref{FTLA}$ covers Lyapunov analysis: first we define finite-time Lyapunov exponents and vectors
(FTLE/Vs) and describe their use for the identification of the tangent
space structure; second we briefly describe the asymptotic theory of partially
hyperbolic sets; third we address the
convergence of the tangent space structure; and fourth we contrast the properties of the FTLE/Vs and their asymptotic
counterparts. In Section $\ref{theory}$ we define
a finite two-timescale set and present the conditions satisfied by
points on a finite-time center manifold. The procedure for applying
the approach is given in Section $\ref{procedure}$. Section $\ref{examples}$ contains detailed examples. Conclusions are given in
Section $\ref{conclusions}$.

%==================================================================

\section{Dynamical System Description and Relevant Geometry} \label{setup}  The methodology we develop will be applied
to
a given coordinate representation of a dynamical system.
Denoting the vector of coordinates by $\mathbf{x}\in \mathbb{R}^n$, in the standard basis with
$2\le n < \infty$, the $\mathbf{x}$-representation of the dynamical
system is
\begin{equation}\label{nldyn} \dot{\mathbf{x}} =
\mathbf{f}(\mathbf{x}) ,\end{equation} where the vector field
$\mathbf{f}:\mathbb{R}^n\rightarrow\mathbb{R}^n$ is a smooth
function. The solution of (\ref{nldyn}) for the initial condition
$\mathbf{x}$ is denoted by $\mathbf{x}(t)=\phi(t,\mathbf{x})$, where
$\mathbf{\phi}(t,\cdot):\mathbb{R}^n\rightarrow\mathbb{R}^n$ is the
$t$-dependent flow associated with the vector field $\mathbf{f}$ and
$\mathbf{\phi}(0,\mathbf{x})=\mathbf{x}$. We assume that
$\mathbf{\phi}$ is complete on $\mathbb{R}^n$ for simplicity, but the methodology developed will only
be applied on a subset of the state space and the properties of the
flow outside this subset are irrelevant.

The linearized dynamics associated with (\ref{nldyn}) are
\begin{equation}\label{lindyn} \dot{\mathbf{v}} = D\mathbf{f}(\mathbf{x})
\mathbf{v} \end{equation} where $D\mathbf{f} := \partial{\mathbf{f}}/\partial{\mathbf{x}}$ and will be analyzed
to characterize the timescales. An initial
point $(\mathbf{x}, \mathbf{v})$ is mapped in time $t$ to the point
$(\mathbf{x}(t),\mathbf{v}(t))=(\phi(t,\mathbf{x}),\Phi(t,\mathbf{x})\mathbf{v})$
where $\Phi$ is the fundamental matrix for the linearized dynamics,
defined such that $\Phi(0,\mathbf{x})=I$, the $n\times n$ identity
matrix. With this initial condition, we refer to $\Phi$ as the transition matrix. Geometrically, for a pair
$(\mathbf{x}, \mathbf{v})$, we view $\mathbf{v}$ as taking values in
the tangent space at $\mathbf{x}$ denoted by
$T_\mathbf{x}\mathbb{R}^n$. The tangent bundle $T\mathbb{R}^n$ is
the union of the tangent spaces over the state space $\mathbb{R}^n$ and
$(\mathbf{x},\mathbf{v})$ is a point in the tangent bundle, with
$\mathbf{v}$ the tangent vector and $\mathbf{x}$ the base point. We
need the interpretation $(\mathbf{x},\mathbf{v})\in T\mathbb{R}^n$,
because the analysis of the linearized dynamics will define a
subspace decomposition of the tangent space and the orientation of
the subspaces will vary with the base point $\mathbf{x}$. Henceforth (\ref{lindyn}) is called the {\em tangent linear
dynamics}.

We adopt the Euclidean metric for $\mathbb{R}^n$ and the Euclidean norm to define the length
of a tangent vector, i.e., for $\mathbf{v}\in
T_\mathbf{x}\mathbb{R}^n$, its length is $\|\mathbf{v}\|=\langle
\mathbf{v},\mathbf{v} \rangle^{1/2}$ and $\langle\cdot,\cdot\rangle$
is the standard inner product.

Let $\mathbf{w}_1, \mathbf{w}_2,\dots,\mathbf{w}_k$, $k\le n$, denote vector
fields, defined on $\mathbb{R}^n$, that vary continuously with
$\mathbf{x}$ and have the property that at each $\mathbf{x}\in
\mathbb{R}^n$, the vectors
$\mathbf{w}_1(\mathbf{x}),\dots,\mathbf{w}_k(\mathbf{x})$ are
linearly independent in $T_\mathbf{x}\mathbb{R}^n$. Then at each
$\mathbf{x}$,  $ \Delta (\mathbf{x})={\rm
span}\{\mathbf{w}_1(\mathbf{x}),\dots,$ $\mathbf{w}_k(\mathbf{x})\}$
is a $k$-dimensional subspace. If $k=n$, then $ \Delta
(\mathbf{x})=T_\mathbf{x}\mathbb{R}^n$ and for each $\mathbf{x}$ the
set of vectors provides a basis for $T_\mathbf{x}\mathbb{R}^n$. If
$k<n$, then $ \Delta (\mathbf{x})$ is a linear subspace of
$T_\mathbf{x}\mathbb{R}^n$; let $ \Delta:=\bigcup_{x\in R^n} \Delta (\mathbf{x})$ denote the subbundle (or distribution)
on $\mathbb{R}^n$. A subbundle is $\Phi$-invariant, if for any
$\mathbf{x}\in \mathbb{R}^n$ and $\mathbf{v}\in
\Delta(\mathbf{x})$, the property $\Phi(t,\mathbf{x})\mathbf{v}\in
\Delta(\phi(t,\mathbf{x}))$ holds for all $t$. Subbundles
$\Delta_1,\dots,\Delta_m$ allow a splitting of the tangent bundle
if $T\mathbb{R}^n=\Delta_1\oplus\dots\oplus
\Delta_m$. If each subbundle in the splitting is
$\Phi$-invariant, then the splitting is a $\Phi$-invariant
splitting.

Let $\mathcal{X}$ be a domain in $\mathbb{R}^n$. A smooth submanifold $\mathcal{M}\subset\mathcal{X}\subset
\mathbb{R}^n$ of dimension
$m<n$ is $\mathcal{X}$-relatively $\phi$-invariant, if for each $\mathbf{x}\in \mathcal{M}$,
$\mathbf{\phi}(t,\mathbf{x})\in \mathcal{M}$ for all $t$ for which $\phi(t,\mathbf{x})$ has not left $\mathcal{X}$. An
equivalent requirement for invariance is that
$\mathbf{f}(\mathbf{x})\in T_\mathbf{x}{\cal M}$ for all
$\mathbf{x}\in {\mathcal{M}}$.

 %==================================================================
 \section{Overview of Approach} \label{overview}

Consider a domain
$\mathcal{X}\subset \mathbb{R}^n$ on which the behavior of (\ref{nldyn}) on a time interval $[0, t_f]$ is of interest.
The tangent linear dynamics (\ref{lindyn}) are analyzed to determine if there is a
splitting of the tangent bundle into stable,
center, and unstable subbundles
$T\mathcal{X}=\mathcal{E}^s\oplus
\mathcal{E}^c\oplus \mathcal{E}^u$ of dimensions $n^s$, $n^c$, and $n^u$, respectively, where the associated exponential
rates indicate that, relative to the time interval $[0, t_f]$, vectors in the stable subbundle $\mathcal{E}^s$ decay
quickly in forward time, vectors in the unstable subbundle $\mathcal{E}^u$ decay quickly in backward time, and vectors in the center subbundle evolve slowly. Then
postulate that there are corresponding invariant manifolds that organize the flow in the state space on the time
interval of interest.
For example, an $n^c$-dimensional invariant center manifold $\mathcal{W}^c\subset\mathcal{X}$ can be postulated. At each
$\mathbf{x}\in\mathcal{W}^c$, $T_\mathbf{x}\mathcal{W}^c = \mathcal{E}^c(\mathbf{x})$
and $\mathbf{f}(\mathbf{x})\in \mathcal{E}^c(\mathbf{x})$.
If
$\{\mathbf{w}_1(\mathbf{x}),\dots,\mathbf{w}_{n-n^c}(\mathbf{x})\}$
is a basis for $[\mathcal{E}^c(\mathbf{x})]^\bot$, the orthogonal complement of $\mathcal{E}^c(\mathbf{x})$, then a
necessary condition for a point $\mathbf{x}\in \mathcal{W}^c$ is the satisfaction of the orthogonality conditions
\begin{equation}\label{orthog}
\langle \mathbf{f}(\mathbf{x}),\mathbf{w}_i(\mathbf{x})\rangle =
0,~~ i=1,...,n-n^c\end{equation} The orthogonality conditions express that $\mathbf{f}(\mathbf{x})$ lies in
$T_\mathbf{x}\mathcal{W}^c$ at each $\mathbf{x}\in \mathcal{W}^c$, i.e., the invariance of $\mathcal{W}^c$.
The orthogonality conditions for $\mathbf{f}$ in (\ref{orthog}) can
be viewed as partial-equilibrium conditions, partial in the sense
that the vector field $\mathbf{f}$ need only be zero when projected
into a certain subspace.  Similarly, orthogonality conditions can be expressed for points on the center-stable
$\mathcal{W}^{cs}$ and center-unstable $\mathcal{W}^{cu}$ manifolds.

\begin{figure}[ht]\centering\includegraphics [scale=.53] {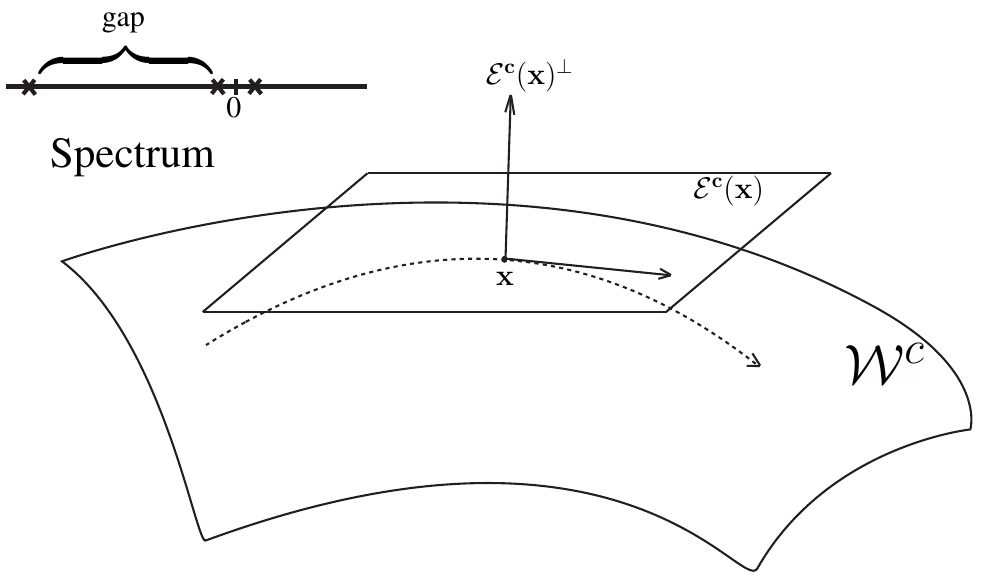}
\caption{Geometry of a two-timescale 3D system with  a 2D normally attracting center manifold. } \label{2Dattracting}
\end{figure}

\begin{figure}[ht]\centering\includegraphics [scale=.53] {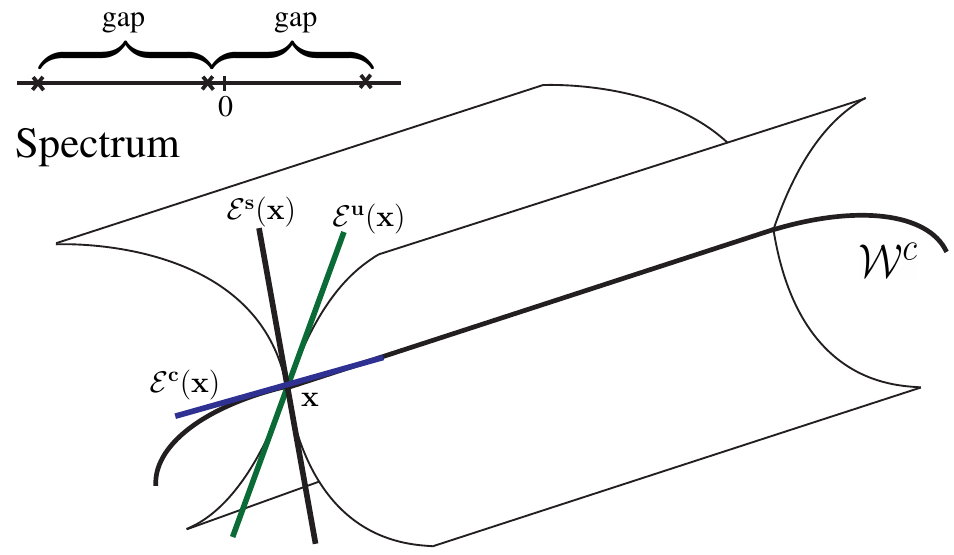}
\caption{Geometry of a two-timescale 3D system with  a 1D normally hyperbolic center manifold.} \label{1Dhyperbolic}
\end{figure}

Figures~\ref{2Dattracting} and \ref{1Dhyperbolic} show examples of center
manifolds in a three-dimensional state space, the relevant geometric objects, and the spectra of characteristic
exponents indicating the exponential rates in the tangent linear dynamics,
consistent with the
geometry. Diagnosing timescale separation and computing such
geometric structure, encompassing both the normally attracting center manifold
(Fig.~\ref{2Dattracting}) and normally hyperbolic center manifold
(Fig.~\ref{1Dhyperbolic}) cases, is our goal. Computationally, determining only low-dimensional manifolds may be feasible, but computing selected points on higher-dimensional manifolds is possible and useful (e.g., \cite{ACC09}).

The theory of partially hyperbolic dynamical systems
\cite{PesinHB,HirschShub} and Oseledec decompositions
\cite{Oseledec} guiding our approach focuses on the behavior on the normally hyperbolic manifold, assumes the manifold is a compact invariant set, and considers behavior on
an infinite-time interval. For many applications, the manifold of interest is not known a priori, so it must be located by analyzing the behavior over a larger region, and the time interval over which the behavior can be analyzed is finite. In this setting, the approach of determining the tangent bundle splitting and
using it to determine the manifold structure in the state space has been
pursued.

A feature that distinguishes the approaches that have been
taken is the means of determining the timescales and splitting. The eigenvalues and
eigenvectors of $D\mathbf{f}(\mathbf{x})$ are used in the ILDM method \cite{MaasPope}. At each point $\mathbf{x}$ of
interest, the eigenvalues provide the timescales and the eigenvectors provide basis vectors for representing the
splitting. When applied to a system in the singularly perturbed normal form, the error in
determining points on a slow manifold using orthogonality conditions formulated with the eigenvectors of
$D\mathbf{f}(\mathbf{x})$ increases as the timescale separation decreases and as
the curvature of the slow manifold increases \cite{Kaper}. The computational singular
perturbation (CSP) method \cite{Lam1,Lam2,KDM96} includes an
iterative procedure that adjusts the eigenvectors of
$D\mathbf{f}(\mathbf{x})$ to basis vectors that better approximate
the slow-fast splitting based on the invariance of these subspaces under
the linear flow; see \cite{KaperCSP} for an error analysis.
Eigen-analysis of the symmetric part of
$D\mathbf{f}(\mathbf{x})$ was employed in \cite{Gorban}, and eigen-analysis of the symmetric part of a reduced form of
$D\mathbf{f}(\mathbf{x})$ characterizing the directions normal to the vector field was used in \cite{Adrover07}. In the
chemical kinetics context when the system is dissipative and
all trajectories asymptotically approach an equilibrium point, a Lyapunov function
is known and a projection to the slow subspace can be derived from it
\cite{Gorban}.

Finite-time Lyapunov analysis is used in \cite{adrover06,JGCD03,AMCpaper}: the FTLEs provide the timescales and the
FTLVs provide the basis vectors for representing the splitting. The FTLEs and the
 FTLVs are the singular values and singular vectors of $\Phi$ for a propagation time $T$.
As the propagation time $T$ goes to zero,
the FTLE/Vs approach the eigenvalues/vectors of the symmetric part
of $D\mathbf{f}(\mathbf{x})$ \cite{Doerner}. For Lyapunov regular points, the limits of the FTLEs, as $T$ goes to
infinity, are the asymptotic
Lyapunov exponents used in the theory of hyperbolic systems \cite{PesinHB,Katok}. Thus FTLA can characterize from instantaneous behavior to average behavior over finite to infinite time intervals, depending on the propagation time used.

In \cite{JGCD03,AMCpaper} FTLA is applied to slow-fast systems to improve the accuracy of slow manifold approximations relative to that of the ILDM method. In \cite{adrover06}, in addition to slow manifolds, FTLA is applied to periodic and chaotic attractors which are outside the domain of applicability of the ILDM and CSP methods when the eigenvectors of $D\mathbf{f}(\mathbf{x})$ either rotate too fast or are complex; see also the earlier work \cite{Lorenz} using FTLA for a chaotic attractor. In the present paper we develop FTLA for partially hyperbolic splittings and normally hyperbolic center manifolds.

For the singularly perturbed normal form depending a small parameter $\varepsilon$, there is an $\varepsilon$-dependent
center manifold $\mathcal{W}^c(\varepsilon)$, which for $\varepsilon=0$ is a manifold composed of equilibrium points
\cite{Fenichel}. For small $\varepsilon$, the flow on $\mathcal{W}^c(\varepsilon)$ is slow. In this case, the nonlinear
dynamics, as well as the tangent linear dynamics, have slow-fast behavior, and it is appropriate to refer to the center
manifold as the {\bf slow manifold}. The ILDM and CSP methods were conceived as $\varepsilon$-free means, as was the FTLA method in \cite{JGCD03,AMCpaper},  of
achieving singular perturbations type results for slow-fast systems; thus it was appropriate to call the center manifold, the slow
manifold. However slow-fast behavior in the tangent linear dynamics does not in general imply that the flow on the
center manifold is slow,\footnote{The authors thank an anonymous reviewer for stimulating our thinking on this issue.} and the approach we develop for normally hyperbolic center manifolds does not require this. Thus we treat a more general class of systems, that includes slow-fast systems, with the common feature that, on the time-interval of interest, there is a fast approach of trajectories to a reduced-order manifold.

%==================================================================
\section{Lyapunov Analysis and Partially Hyperbolic Sets - Finite-Time versus Asymptotic} \label{FTLA}
 In this section we present the methodology for
characterizing the tangent linear dynamics (\ref{lindyn}), along
trajectories of the nonlinear system (\ref{nldyn}), that will be used to define and diagnose two-timescale behavior. We
refer to
this methodology as {\em Lyapunov analysis}. Because we need to determine, in a limited finite-time, a good approximation of an invariant splitting that in principle requires asymptotic Lyapunov analysis, we need to define the finite-time splitting in a way that will converge as fast as possible towards the desired invariant splitting. We clarify that defining the splitting in terms of FTLVs accomplishes this. In the
first subsection, we present a finite-time version of Lyapunov
analysis, modeled after the asymptotic version described in Barreira
and Pesin \cite{Pesin} and Katok and Hasselblatt \cite{Katok}.
In the second subsection, we describe how asymptotic Lyapunov exponents or vectors can be used to define the ideal invariant splittings; in the
third subsection, the convergence rate of a Lyapunov subspace is
characterized; and in the final subsection, the products of asymptotic and finite-time Lyapunov analysis are contrasted---in
preparation for the finite-time approach
presented in the remaining sections. See also \cite{adrover06,Legras,JGCD03,Samelson,Yoden} for presentations of
finite-time Lyapunov analysis.

\subsection{Finite-Time Lyapunov Exponents/Vectors and Tangent Space Structure}\label{FTLE}

The forward and backward FTLEs for a vector $\mathbf{v}\in T_{\mathbf{x}}\mathbb{R}^n$ are given by
\begin{equation}\label{forwardexp}
\mu^{\pm}(T,\mathbf{x}, \mathbf{v}) := \frac{1}{T}\ln\Lambda^{\pm}(T,\mathbf{x}, \mathbf{v})=
\frac{1}{T}\ln\frac{\|\Phi(\pm T,\mathbf{x}) \mathbf{v}\|}{\|\mathbf{v}\|} \\
\end{equation}
where $T$ is the propagation time, also referred to as the averaging time, and is always taken to be
positive whether the propagation is forward or backward.
 Variables computed by
forward and backward propagation are labeled with
superscripts $^+$ and $^-$ respectively.
For $\mathbf{v}=\mathbf{0}$, define
$\mu^+(T,\mathbf{x},\mathbf{0})$ $ =\mu^-(T,\mathbf{x},\mathbf{0})$
$=-\infty$. The FTLE is the average exponential rate of growth/decay over the time interval $[0,T]$.

Discrete forward and backward Lyapunov spectra, for each
$(T,\mathbf{x})$, can be defined as follows. Define
$\mathbf{l}^+_i(T,\mathbf{x})$, $i=1,\dots,n$, to be an orthonormal basis of
$T_\mathbf{x}\mathbb{R}^n$ with the minimum sum of exponents, i.e.,
the minimum value of $\Sigma_{i=1}^{n}
\mu^+_i(T,\mathbf{x},\mathbf{l}^+_i(T,\mathbf{x}))$ over all
orthonormal bases \cite{Dieci}. The forward Lyapunov spectrum is the
set of exponents corresponding to the minimizing solution, namely,
$\{\mu^+_i(T,\mathbf{x}), i=1,\dots,n\}$. The FTLEs are assumed to be in ascending order. The Lyapunov spectrum is
unique, though the minimizing basis is not in general.

Geometrically, the unit
$n$-sphere centered at the origin in $T_\mathbf{x}\mathbb{R}^n$
propagates under the tangent linear dynamics to an $n$-dimensional
ellipsoid in $T_{\phi(T,\mathbf{x})}\mathbb{R}^n$; the principal
semi-axes of the ellipsoid are
$\exp[\mu_i^+(T,\mathbf{x})T]\mathbf{n}^+_i(T,\mathbf{x})$,
$i=1,\dots,n$ and the unit vectors in $T_\mathbf{x}\mathbb{R}^n$
that evolve to these vectors are respectively
$\mathbf{l}^+_i(T,\mathbf{x}),i=1,\dots,n$.

The backward Lyapunov spectrum $\{\mu^-_i, i=1,\dots,n\}$ consists of the exponents
for the unit vectors $\mathbf{l}^-_i(T,\mathbf{x}), i=1,\dots,n$ in $T_\mathbf{x}\mathbb{R}^n$ that map to
principal axes of an $n$-ellipsoid in
$T_{\phi(-T,\mathbf{x})}\mathbb{R}^n$. Descending order is assumed for the backward FTLEs.

\begin{figure}[ht]\centering\includegraphics[scale=.46] {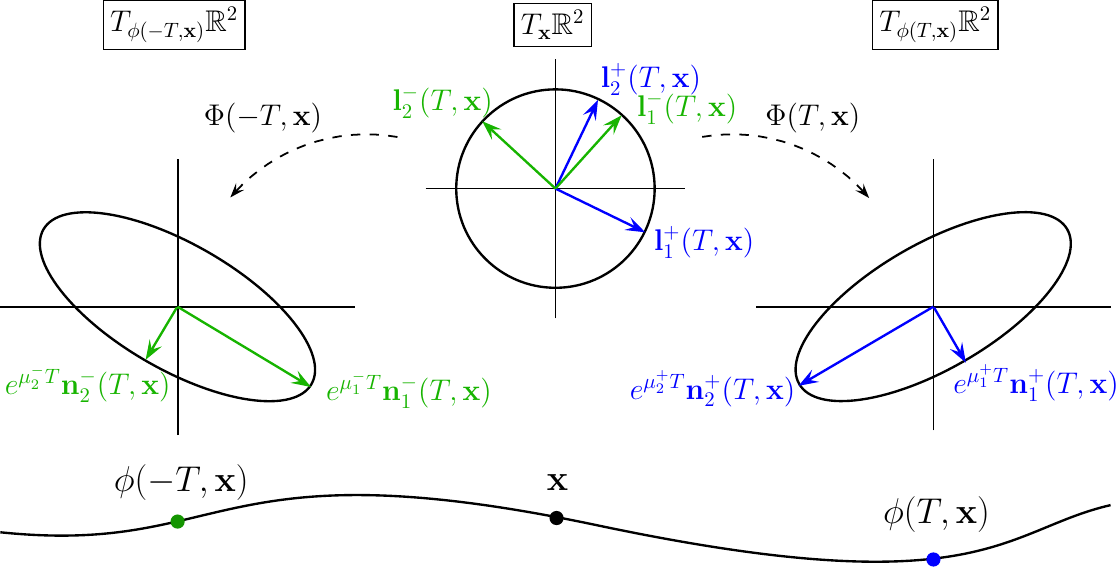}
\caption{Trajectory of nonlinear system and associated tangent
spaces, illustrating the role of the Lyapunov exponents and vectors
in the forward and backward propagation of a sphere of tangent
vectors. Blue objects correspond to forward propagation, and green
objects correspond to backward propagation. The arguments $(T,
\mathbf{x})$ of the FTLE/Vs have been suppressed.}
\label{ellipsoids1}
\end{figure}

 The
$\mathbf{l}_i^+(T,\mathbf{x})$ and the
$\mathbf{l}_i^-(T,\mathbf{x})$ vectors, for $i=1,\dots,n$, referred
to as forward and backward FTLVs, respectively, will be used to
define subspaces in $T_\mathbf{x}\mathbb{R}^n$ associated with
different exponential rates.See Fig.~\ref{ellipsoids1} for the case of $n=2$.

 We assume that the FTLEs are always distinct, i.e., non-degenerate. This assumption simplifies the presentation and is
needed in slightly stronger form for
the
 subspace convergence proof presented in \ref{appendixA}. We
 note that distinctness is also related to integral
 separation and the stability of the Lyapunov exponents with respect
 to perturbations in the linearized system matrix, $D\mathbf{f}(\mathbf{x})$
 \cite{Dieci}. Later we accommodate degeneracies in an initial ``transient" phase that is short relative to
 the time interval under consideration by modifying the assumption to
 hold for $T\ge t_s$, for an appropriate value of the time $t_s$.

 The
following subspaces, for $i=1,\dots,n$, can be defined by the
orthonormal FTLVs
\begin{equation}\label{Ldefs}\begin{array}{ccl} {\cal L}^+_i(T,\mathbf{x}) &:=&
span\{\mathbf{l}^+_1(T,\mathbf{x}),
\dots,\mathbf{l}^+_{i}(T,\mathbf{x})\} ,\\ {\cal
L}^-_i(T,\mathbf{x}) &:=& span\{\mathbf{l}^-_i(T,\mathbf{x}),
\dots,\mathbf{l}^-_{n}(T,\mathbf{x})\} , \end{array}\end{equation}
and will be referred to as finite-time Lyapunov subspaces. For any
{$i \in \{1,2,\dots,n\}$},
 $\mu^{{+}}(T,\mathbf{x},\mathbf{v})\le
\mu^{{+}}_i(T,\mathbf{x})$ for any $\mathbf{v}\in
 {\cal L}^{{+}}_i(T,\mathbf{x})$. However, for finite $T$, there also exist vectors $\mathbf{v}\in
 T_\mathbf{x}\mathbb{R}^n\setminus {\cal L}^{{+}}_{i}{(T,\mathbf{x})}$ for which $\mu^{{+}}(T,\mathbf{x},\mathbf{v})\le
\mu^{{+}}_i(T,\mathbf{x})$. Analogous properties hold for the backward-time exponents and subspaces.

If a collection of $r\le n$ linear subspaces of
$T_\mathbf{x}\mathbb{R}^n$ can be ordered such that $
\Delta_{1}(\mathbf{x})\subset
\Delta_{2}(\mathbf{x})\subset\dots\subset
\Delta_r(\mathbf{x})=T_\mathbf{x}\mathbb{R}^n$ with all inclusions
strict, then this collection of nested subspaces defines a {\it
filtration} of $T_\mathbf{x}\mathbb{R}^n$. The nested sequences of
subspaces
 \begin{equation}\label{forwardfiltration}
\{\mathbf{0}\}=: {\cal L}_0\subset
 {\cal L}_{1}^+(T,\mathbf{x})\subset {\cal L}_{2}^+(T,\mathbf{x})\subset\dots\subset
 {\cal L}_n^+(T,\mathbf{x})=T_{\mathbf{x}}\mathbb{R}^n ,
\end{equation}
  \begin{equation}\label{backwardfiltration}
T_\mathbf{x}\mathbb{R}^n=  {\cal L}^-_{1}(T,\mathbf{x})\supset {\cal
L}^-_{2}(T,\mathbf{x})\supset\dots\supset
 {\cal L}^-_{n}(T,\mathbf{x})\supset  {\cal
 L}^-_{n+1}:=\{\mathbf{0}\} ,
\end{equation} are forward and backward filtrations \cite{Pesin,Katok} of $T_{\mathbf{x}}\mathbb{R}^n$.

 We need both forward and backward filtrations, because their intersections
 are of particular interest, as motivated by the following.
 Consider a two-dimensional nonlinear system with  an equilibrium point $\mathbf{x}_{e}$.
Assume the linearized dynamics at $\mathbf{x}_{e}$ are characterized
by distinct eigenvalues $\lambda_1$ and $\lambda_2$, with
$\lambda_1<\lambda_2<0$, and corresponding unit eigenvectors
$\mathbf{e}_1$ and $\mathbf{e}_2$. As $T\rightarrow \infty$, the
FTLEs at $\mathbf{x}_{e}$ approach the eigenvalues, i.e., $\mu^+_1(T,\mathbf{x}_{e})\rightarrow \lambda_1$
and $\mu^+_2(T,\mathbf{x}_{e})\rightarrow \lambda_2$, and the first Lyapunov vector
approaches the corresponding eigenvector $\mathbf{l}^+_1(T,\mathbf{x}_{e})\rightarrow
\mathbf{e}_1$. The second Lyapunov vector $\mathbf{l}^+_2(T,\mathbf{x}_{e})$
approaches $\mathbf{e}_1^\bot$, the vector perpendicular to
$\mathbf{e}_1$. The subspace ${\cal L}^+_1(T,\mathbf{x}_{e})$ thus
approaches $\mathcal{E}^1(\mathbf{x}_{e})=span\{\mathbf{e}_1\}$, the
eigenspace for $\lambda_1$ as $T\rightarrow \infty$, whereas ${\cal
L}^+_2(T,\mathbf{x}_{e})=T_{\mathbf{x}_{e}}\mathbb{R}^2$ for any
$T$. It is desired instead to obtain the invariant splitting
$T_{\mathbf{x}_{e}}\mathbb{R}^2=\mathcal{E}^1(\mathbf{x}_{e})\oplus
\mathcal{E}^2(\mathbf{x}_{e})$ where
$\mathcal{E}^2(\mathbf{x}_{e})=span\{\mathbf{e}_2\}$. However, asymptotically
all the vectors not in ${\cal L}^+_1$ have the Lyapunov exponent
$\mu^+_2=\lambda_2$; thus the Lyapunov exponents for forward-time
propagation do not distinguish $\mathcal{E}^2$.
 The way to obtain $\mathcal{E}^2$ is by
repeating the same analysis for backward-time propagation; in this
case, the situation is reversed: asymptotically
$\mathbf{l}^-_2(T,\mathbf{x}_{e})\rightarrow \mathbf{e}_2$ and $\mathcal{E}^2$ can be
distinguished, whereas $\mathcal{E}^1$ cannot \cite{Katok,Young}.

  %=================================================================
  \subsection{Asymptotic Lyapunov Analysis and Partially Hyperbolic Set}\label{asympt}

  We draw from \cite{Pesin,PesinHB} to present the asymptotic theory,
  covering only those definitions and results that serve to motivate
  and support
  our definitions and results for the finite-time case. Asymptotic Lyapunov analysis was introduced in
\cite{Lyapunov} and related to tangent space geometry in
\cite{Oseledec}. The theory of partially hyperbolic sets is described in \cite{PesinHB} where references to the original
work are given. The definition of a uniform partially hyperbolic set given next requires exponential bounds uniformly,
i.e., on all time intervals for a given trajectory as well as for all trajectories in the set.

\begin{definition}\label{PHSdef} \cite{PesinHB} A compact $\phi$-invariant set $\mathcal{Y}\subset\mathbb{R}^n$ is a
{\bf  uniform partially
hyperbolic set}, if there exists a $\Phi$-invariant
splitting
  \begin{equation}\label{splitting} T_\mathbf{x}{\mathbb{R}^n}=\mathcal{E}^s(\mathbf{x})\oplus
  \mathcal{E}^c(\mathbf{x})\oplus \mathcal{E}^u(\mathbf{x}) \end{equation}
 on $\mathcal{Y}$ and numbers $\sigma$, $\nu$, and $C$, with $0<\sigma<\nu$ and
 $1\le C<\infty$,
 such that $\forall t>0$
 \begin{equation}\label{expbdds}
 \begin{array}{lcl} \mathbf{v}\in
\mathcal{E}^s(\mathbf{x}) &\Rightarrow &
\|\Phi(t,\mathbf{x})\mathbf{v}\|\le C e^{-\nu t}\|\mathbf{v}\|,
\\
\mathbf{v}\in \mathcal{E}^c(\mathbf{x}) &\Rightarrow & C^{-1}e^{-\sigma
t}\|\mathbf{v}\|\le \|\Phi(t,\mathbf{x})\mathbf{v}\|\le C e^{\sigma t}\|\mathbf{v}\|,\\
 \mathbf{v}\in \mathcal{E}^u(\mathbf{x}) &\Rightarrow &
\|\Phi(-t,\mathbf{x})\mathbf{v}\|\le C e^{-\nu t}\|\mathbf{v}\|.
\end{array} \end{equation} \end{definition}

Consistent with the definition, consider for the moment a compact, invariant set $\mathcal{Y}\subset \mathbb{R}^n$.
When the infinite-time limits $(T \rightarrow \infty)$ of the
exponents in (\ref{forwardexp}) exist at $\mathbf{x}\in \mathcal{Y}$
for all $\mathbf{v}\in T_\mathbf{x}\mathbb{R}^n$, they are denoted by
$\mu^+(\mathbf{x},\mathbf{v})$ and $\mu^-(\mathbf{x},\mathbf{v})$
and the system is said to be, respectively, \emph{forward regular}
and \emph{backward regular} at $\mathbf{x}$. There are at most $n$
distinct exponents for the vectors in
$T_\mathbf{x}{\mathbb{R}^n}\backslash\{0\}$. Consistent with our assumption
for the finite-time case, assume that there are $n$ distinct
exponents, denoted $\mu^+_i(\mathbf{x})$, $i=1,\dots,n$ for
forward time and $\mu^-_i(\mathbf{x})$, $i=1,\dots,n$ for
backward time, with the forward exponents in ascending order and the
backward exponents in descending order. Lyapunov subspaces are
defined by $\mathcal{L}^+_i(\mathbf{x}):=\{\mathbf{v}\in
T_\mathbf{x}\mathbb{R}^n:
\mu^+(\mathbf{x},\mathbf{v})\le\mu^+_i(\mathbf{x})\}$ and
$\mathcal{L}^-_i(\mathbf{x}):=\{\mathbf{v}\in T_\mathbf{x}\mathbb{R}^n:
\mu^-(\mathbf{x},\mathbf{v})\le\mu^-_i(\mathbf{x})\}$.  Forward and
backward filtrations are defined as in (\ref{forwardfiltration}) and
(\ref{backwardfiltration}) using the asymptotic Lyapunov subspaces. The system is \emph{Lyapunov regular}
\cite{Pesin} at $\mathbf{x}$ if (i) it is forward and backward
regular at $\mathbf{x}$, (ii)
$\mu^+_i(\mathbf{x})=-\mu^-_i(\mathbf{x})$, $i=1,\dots,n$, (iii) the
forward and backward filtrations have the same dimensions, (iv)
there exists a splitting
$T_\mathbf{x}\mathcal{Y}=\mathcal{E}^1(\mathbf{x})\oplus\dots\oplus \mathcal{E}^n(\mathbf{x})$
into invariant subspaces such that
$\mathcal{L}^+_i(\mathbf{x})=\mathcal{E}^1(\mathbf{x})\oplus\dots\oplus
\mathcal{E}^i(\mathbf{x})$ and
$\mathcal{L}^-_i(\mathbf{x})=\mathcal{E}^i(\mathbf{x})\oplus\dots\oplus
\mathcal{E}^n(\mathbf{\mathbf{x}})$, $i=1,\dots,n$, and (v) for any
$\mathbf{v}\in \mathcal{E}^i(\mathbf{x})\setminus\{0\}$,
$\lim_{t\rightarrow\pm\infty}(1/|t|)\ln\|\Phi(t,\mathbf{x})\mathbf{v}\|=\mu^{\pm}_i(\mathbf{x})$.
The invariant splitting described in (iv) and (v) is referred
to as Oseledec's decomposition.

Next we describe how the Lyapunov exponents and vectors can be used to diagnose and specify a uniform partially hyperbolic set. For the
purpose of motivating the finite-time theory presented in the next
section, assume the system (\ref{nldyn}) is Lyapunov regular at
all the points of a compact, invariant set $\mathcal{Y}$. Suppose we find that at each $\mathbf{x}\in
\mathcal{Y}$, there are, $n^s$ large negative exponents, $n^c$
small in absolute value exponents, and $n^u$ large positive
exponents, with $n^s+ n^c+n^u=n$. That is,
uniformly in $\mathbf{x}$, there is a splitting of the forward Lyapunov spectrum $sp^+(\mathbf{x})$ of the form
$sp^+(\mathbf{x}):=sp^s(\mathbf{x})\cup sp^c(\mathbf{x})\cup
sp^u(\mathbf{x})$ where
$sp^s(\mathbf{x}):=\{\mu^+_1(\mathbf{x}),\dots,\mu^+_{n^s}(\mathbf{x})\}$,
$sp^c(\mathbf{x}):=\{\mu^+_{n^s+1}(\mathbf{x}),\dots,\mu^+_{n^s+n^c}(\mathbf{x})\}$,
and
$sp^u(\mathbf{x}):=\{\mu^+_{n^s+n^c+1}(\mathbf{x}),\dots,\mu^+_{n}(\mathbf{x})\}$.
We can construct a
$\Phi$-invariant splitting with
\begin{equation}\label{Edefs}
\begin{array}{lcl}  \mathcal{E}^s(\mathbf{x})&=&
{\cal L}^+_{n^s}(\mathbf{x}),\\  \mathcal{E}^c(\mathbf{x})&=&
{\cal L}^+_{n^s+n^c}(\mathbf{x})\cap {\cal L}^-_{n^s+1}(\mathbf{x}),\\
 \mathcal{E}^u(\mathbf{x})&=& {\cal L}^-_{n^s+n^c+1}(\mathbf{x}).
\end{array} \end{equation}
Although Lyapunov vectors are not normally used to define the subspaces in
the asymptotic theory, they can be as follows. Let $\{\mathbf{l}^+_i(\mathbf{x}),
i=1,\dots,n\}$ denote an orthonormal basis for $T_\mathbf{x}\mathbb{R}^n$
such that $\{\mathbf{l}^+_j(\mathbf{x}), j=1,\dots,i\}$ is a basis
for $\mathcal{L}^+_i(\mathbf{x})$ for
 $ i=1,\dots, n$. Let
$\{\mathbf{l}^-_i(\mathbf{x}), i=1,\dots,n\}$ denote an orthonormal
basis for $T_\mathbf{x}\mathbb{R}^n$ such that $\{\mathbf{l}^-_j(\mathbf{x}),
j=i,\dots,n\}$ is a basis for $\mathcal{L}^-_i(\mathbf{x})$ for
 $ i=1,\dots, n$. When there are $n$ distinct Lyapunov exponents
as we are assuming, it follows that these bases are unique up to
multiplication of individual vectors by $\pm 1$. These are clearly
the orthonormal bases that minimize the sum of the asymptotic
exponents over the set of orthonormal bases, and hence the basis
vectors are the asymptotic counterparts of the FTLVs.

The final step in specifying the uniform partially hyperbolic set is to define the constants $\sigma=\sigma_0 +
\varepsilon$ and $\nu=\nu_0 - \varepsilon$ where $\varepsilon>0$ is an arbitrarily small constant,
\begin{equation}\begin{array}c
\sigma_0 :=\max\{|\overline{\mu}^c|,|\underline{\mu}^c|\}, \nu_0 :=\min\{-\overline{\mu}^s,\underline{\mu}^u\},
\end{array}\end{equation}\label{mus}
and
\begin{equation}
\begin{aligned}
\overline{\mu}^s=\sup\limits_{\mathbf{x}\in \mathcal{Y}} \mu^+_{n^s}(\mathbf{x}), \hspace{14mm} &
\overline{\mu}^c=\sup\limits_{\mathbf{x}\in \mathcal{Y}} \mu^+_{n^s+n^c}(\mathbf{x}),\\
\underline{\mu}^u=\inf\limits_{\mathbf{x}\in \mathcal{Y}} \mu^+_{n^s+n^c+1}(\mathbf{x}),\hspace{5mm} &
\underline{\mu}^c=\inf\limits_{\mathbf{x}\in \mathcal{Y}} \mu^+_{n^s+1}(\mathbf{x}).
\end{aligned}
\end{equation}

\noindent The bounds are specified in terms of the forward-time exponents $\mu^+$ as defined in (\ref{forwardexp}), but
given the
property (ii) of Lyapunov regularity, the backward-time exponents could have been used. For a partially hyperbolic set
we must have $0 <
\sigma_0 < \nu_0$. Then for sufficiently small $\varepsilon$, there exists a positive, finite
constant $C$ such that the bounds (\ref{expbdds}) hold.

\subsection{Exponential Lyapunov Subspace Convergence} In this
subsection, we relate the finite-time tangent space structure
introduced in Section \ref{FTLE} to the asymptotic tangent space structure
described in Section \ref{asympt}.  Basically the important subspaces converge exponentially fast to their asymptotic
counterparts, and it is this property that makes FTLA viable.

We need to consider the
distance between the subspaces ${\cal L}_j^+(T_1,\mathbf{x})$ and
${\cal L}_j^+(T_2,\mathbf{x})$ in $T_\mathbf{x}\mathbb{R}^n$. For any value of
$j$ in the index set $\{1,2,\dots,n\}$,
let $L_j^+(T,\mathbf{x})$ denote the
matrix whose columns are the Lyapunov vectors
$\mathbf{l}_i^+(T,\mathbf{x})$,~$i=1,\dots,j$, and
$L_{j'}^+(T,\mathbf{x})$ denote the matrix whose columns are the
Lyapunov vectors $\mathbf{l}_i^+(T,\mathbf{x})$,~$i=j+1,\dots,n$.
Then the distance between the subspaces ${\cal
L}_j^+(T_1,\mathbf{x})$ and ${\cal L}_j^+(T_2,\mathbf{x})$ is
\begin{equation}\label{distance}\begin{array}{ll}\textrm{dist}({\cal L}_j^+(T_1,\mathbf{x}),
{\cal L}_j^+(T_2,\mathbf{x}))&=\|L_j^+(T_1,\mathbf{x})^T
L_{j'}^+(T_2,\mathbf{x})\|_2\\&=\|L_{j}^+(T_2,\mathbf{x})^T
L_{j'}^+(T_1,\mathbf{x})\|_2.\end{array}\end{equation}
This result is a special case of Theorem
2.6.1 in \cite{golub}, page 76, and the facts that the columns of
$L_j^+(T,\mathbf{x})$ provide an orthogonal basis for ${\cal
L}_j^+(T,\mathbf{x})$ and the columns of ${L}_{j'}^+(T,\mathbf{x})$
are mutually orthogonal to the columns of ${L}_j^+(T,\mathbf{x})$.

At a forward regular point for which there
exists $t_s>0$ such that for $T>t_s$ there is a
nonzero lower bound ${\Delta\mu}^+_j(\mathbf{x})$ on the spectral gap
$\mu^+_{j+1}(T,\mathbf{x})-\mu^+_{j}(T,\mathbf{x})$, for a
specific value of $j$, the subspace $ {\cal L}^+_j(T,\mathbf{x})$
approaches the fixed subspace ${\cal L}^+_j(\mathbf{x})$, defined in
Section \ref{asympt} in terms of the asymptotic Lyapunov exponent
$\mu^+_j(\mathbf{x})$. It approaches it at an exponential rate characterized, for
every sufficiently small $\Delta T
>0$, by
\begin{equation}\label{expbdd} \textrm{dist}({\cal L}_j^+(T,\mathbf{x}),
{\cal L}_j^+(T+\Delta T,\mathbf{x}))\le K
e^{-{\Delta\mu}^+_j(\mathbf{x})\cdot T},\end{equation} for
all $T>t_s$, where $K>0$ is $\Delta T$ dependent but $T$
independent. Similarly, as $T$ increases, the subspace $ {\cal
L}^-_k(T,\mathbf{x})$ approaches the fixed subspace ${\cal
L}^-_k(\mathbf{x})$ at a rate proportional to
$exp(-\Delta\mu^-_k(\mathbf{x}) \cdot T)$ where $\Delta\mu^-_k(\mathbf{x})$ is the spectral gap lower bound for backward
propagation and $T > t_s$.
For the technical details see \ref{appendixA}.

\subsection{Differences Between Finite-Time and Asymptotic Lyapunov Analysis}
As discussed in Section \ref{asympt}, in the asymptotic setting either Lyapunov exponents or vectors
can serve to define the Lyapunov subspaces and
tangent space splitting, and the results are equivalent. In
contrast, the FTLEs and FTLVs define different tangent space
objects. For example, if one defines the $i^{th}$ forward finite-time Lyapunov
subspace at $\mathbf{x}$ as $\mathcal{V}_i^+(T,\mathbf{x}):=\{\mathbf{v}\in
T_\mathbf{x}\mathbb{R}^n:
\mu^+(T,\mathbf{x},\mathbf{v})\le\mu_i^+(T,\mathbf{x})\}$, one gets
not a subspace, but an object with non-zero volume centered on the FTLV-defined Lyapunov
subspace $\mathcal{L}_i^+(T,\mathbf{x})$ (also noted in \cite{Wiggins}). To see this, consider the
tangent vector $\mathbf{v}=\mathbf{u} + \beta \mathbf{w}$ in
$T_\mathbf{x}\mathbb{R}^n$, with $\mathbf{u}\in
\mathcal{L}_i^+(T,\mathbf{x})$, $\mathbf{w}\in
[\mathcal{L}_i^+(T,\mathbf{x})]^\bot$, and $\beta$ a scalar
constant. For a given $T$, there exist nonzero values of $\beta$
close enough to zero that $\mathbf{v}$ will belong to
$\mathcal{V}_i^+(T,\mathbf{x})$, whereas it does not belong to
$\mathcal{L}_i^+(T,\mathbf{x})$. Under certain conditions \cite{Oseledec}, as $T$
increases, $\mathcal{L}_i^+(T,\mathbf{x})$ converges to its
asymptotic value $\mathcal{L}_i^+(\mathbf{x})$  and $\mathcal{V}_i^+(T,\mathbf{x})$ converges to
$\mathcal{L}_i^+(T,\mathbf{x})$ and thus to
$\mathcal{L}_i^+(\mathbf{x})$ as well. Because the FTLV-defined Lyapunov subspace convergence is exponential in $T$,
while the Lyapunov
exponent convergence is much slower, perhaps proportional to $1/T$
\cite{Orszag}, in the finite-time setting we define the Lyapunov subspaces in terms of the FTLVs.

 The asymptotic Lyapunov exponents for Lyapunov regular points exist as limits, are metric
independent, are constant on a trajectory, and include a zero exponent associated with the vector field direction.
These properties are not shared in general by the FTLEs. The FTLEs
depend on $\mathbf{x}$ and $T$; there need not be a zero exponent associated with the vector field direction. FTLEs can indicate local behavior which, if not uniformly present, would not be indicated by the asymptotic Lyapunov exponents. Another potential
feature in the FTLEs is ``nonmodal behavior" \cite{Schmid} which has required the introduction of the delayed start time $t_s\ge
0$ to avoid a brief initial transient, relative to the time interval of interest, during which  the FTLEs can be quite different than they will be for even moderate finite times. FTLEs
are in general  metric dependent.  In the
present paper, we use the Euclidean metric exclusively, though any
Riemmanian metric could be used  \cite{GK2,LekienRoss,JGCD03}. If finite-time two-timescale  behavior is not
present in the original metric under consideration, there may be
another metric for which there is two-timescale behavior, as noted
by Greene and Kim \cite{GK2}.

%=================================================================
\section{Finite-Time Two-Timescale Set and Center Manifold - Theory} \label{theory}
We identify the potential for manifold structure in a state-space region by determining if a representative set ${\mathcal{X}}\subset \mathbb{R}^n$ is a finite-time uniform two-timescale set. A
two-timescale set has a special tangent space structure and allows
us to formulate invariance-based orthogonality conditions that would be satisfied at points of center, center-stable, and center-unstable manifolds, if such manifolds are present. For the purpose of defining and diagnosing
two-timescale behavior, $\mathcal{X}$ could be a point or a segment of a trajectory, as examples, but in
the search for manifold structure, $\mathcal{X}$ is typically a domain of the state space. The domain is typically not $\phi$-invariant, so it is crucial to clarify what information is required and how much time it takes to resolve it. And because only limited integration time is available, the definition of a finite-time two-timescale set must account for finite-time features that are of no consequence in asymptotic Lyapunov theory.

\subsection{Finite-Time Two-Timescale Set} \label{2TS set}
Definition \ref{definition1} of a  finite-time uniform two-timescale set is guided by Def. \ref{PHSdef} of a
uniformly partially hyperbolic set and consideration of convergence. Several time constants\footnote{For an exponential function of time,
$e^{\kappa
t}$, the time constant $|\kappa|^{-1}$ is the time $t$ at which
the function equals $e^{+1}$ or $e^{-1}$ as appropriate for the sign of $\kappa$.} play key roles. The spectral gap
$\Delta\mu$ must be large enough relative to the common available maximum averaging time $\overline{T}$ that the tangent
space
splitting can be accurately resolved; hence the convergence time constant $\Delta\mu^{-1}$ should allow the
finite-time subspaces to
converge over at least several time constants toward their ideal
infinite-time limits. The fast and slow time constants (i.e., timescales), $\nu^{-1}$ and $\sigma^{-1}$, appear in the
bounds that characterize the disparate exponential rates in the tangent linear dynamics, as further interpreted in
Section \ref{Sig}.

{\begin{definition}\label{definition1}
A set $\mathcal{X}\subset
\mathbb{R}^n$, $n\ge 2$, is a \textbf{uniform finite-time two-timescale set} for
(\ref{nldyn}) with respect to the Euclidean metric, with fast time constant $\nu^{-1}$ and slow time constant
$\sigma^{-1}$, and convergence time constant $\Delta\mu^{-1}$, resolvable over
 $\Delta\mu (\overline{T}-t_s)$ convergence time constants, if there exist positive integers  $n^s$, $n^c$
and $n^u$, with $n^s+n^c+n^u=n$, a delayed start time $t_s$, a cut-off time $t_c$, and an available averaging time
$\overline{T}$ with $0\le t_s< t_c\le\overline{T}$ such that the following three properties are
satisfied. We use the notation $\mathcal{T}=(t_s,\overline{T}]$ and
$\mathcal{T}_c=(t_s,t_c]$.
\begin{enumerate}
\item Uniform Spectral Gaps --
There exist positive constants $\alpha$ and $\beta$ with $\beta-\alpha$ $>0$ such that, uniformly on $\mathcal{T}\times
\mathcal{X}$, the forward and backward Lyapunov spectra are separated by gaps of size $\Delta\mu=\beta-\alpha$  into
$n^{s}$, $n^c$ and $n^{u}$ dimensional subsets as illustrated in Fig.~\ref{slow_spectra_fig} and specified by
\begin{equation}
\begin{aligned}
 \mu^+_{n^{s}}\le - \beta,&\hspace{2mm}      -\alpha\le\mu^+_{n^{s}+1},\hspace{2mm}
 \mu^+_{n^{s}+n^{c}}\le\alpha,\hspace{2mm}    \beta\le\mu^+_{n^{s}+n^{c}+1},\\
-\mu^-_{n^{s}}\le-\beta,&\hspace{2mm}    -\alpha\le-\mu^-_{n^{s}+1},\hspace{2mm}
-\mu^-_{n^{s}+n^{c}}\le\alpha,\hspace{2mm}   \beta\le-\mu^-_{n^{s}+n^{c}+1}.
\end{aligned}
\end{equation}
\item Tangent Bundle Splitting -- On $\mathcal{X}$, there is a continuous splitting \begin{equation}\label{FTsplitting}
T_\mathbf{x}{\mathbb{R}^n}=\mathcal{E}^s(\overline{T},\mathbf{x})\oplus
  \mathcal{E}^c(\overline{T},\mathbf{x})\oplus \mathcal{E}^u(\overline{T},\mathbf{x}) ,\end{equation}
 where \begin{equation}\label{EdefsT}
\begin{array}{lcl}  \mathcal{E}^s(\overline{T},\mathbf{x})&=&
{\cal L}^+_{n^s}(\overline{T},\mathbf{x}),\\  \mathcal{E}^c(\overline{T},\mathbf{x})&=&
{\cal L}^+_{n^s+n^c}(\overline{T},\mathbf{x})\cap {\cal L}^-_{n^s+1}(\overline{T},\mathbf{x}),\\
 \mathcal{E}^u(\overline{T},\mathbf{x})&=& {\cal L}^-_{n^s+n^c+1}(\overline{T},\mathbf{x}).
\end{array} \end{equation}
\item Two Timescales -- There exist positive numbers $\nu$ and $\sigma$ with $\nu > \sigma$ such that at each
$\mathbf{x}\in \mathcal{X}$
for all $t\in \mathcal{T}_c$
\begin{equation}
\begin{aligned}
\mathbf{v}\in \mathcal{E}^s(\overline{T},\mathbf{x})\Rightarrow &
\begin{cases}
\|\Phi(-t,\mathbf{x})\mathbf{v}\| \ge e^{\nu t} \|\mathbf{v}\| & \\
\|\Phi(t,\mathbf{x})\mathbf{v}\| \le e^{-\nu t} \|\mathbf{v}\| &
\end{cases},\\
\mathbf{v}\in \mathcal{E}^c(\overline{T},\mathbf{x})\Rightarrow &
\begin{cases}
e^{-\sigma t} \|\mathbf{v}\|\le\|\Phi(t,\mathbf{x})\mathbf{v}\| \le e^{\sigma t} \|\mathbf{v}\| & \\
e^{-\sigma t} \|\mathbf{v}\|\le\|\Phi(- t,\mathbf{x})\mathbf{v}\| \le e^{\sigma t} \|\mathbf{v}\| &
\end{cases},\\
\mathbf{v}\in \mathcal{E}^u(\overline{T},\mathbf{x})\Rightarrow &
\begin{cases}
\|\Phi(-t,\mathbf{x})\mathbf{v}\| \le e^{-\nu t}\|\mathbf{v}\| & \\
\|\Phi(t,\mathbf{x})\mathbf{v}\| \ge e^{\nu t}\|\mathbf{v}\| &
\end{cases}.\\
\end{aligned}
\end{equation}
\end{enumerate}
It is assumed that $n^c\ge 1$. Either $n^s$ or $n^u$ can be zero, but not both.
 For $n^s=0$,
   $\mathcal{E}^s$
 is not relevant; similarly, for $n^u=0$,
 $\mathcal{E}^u$ is not relevant.
 \end{definition}

\begin{figure}[ht]\centering\includegraphics [scale=.642] {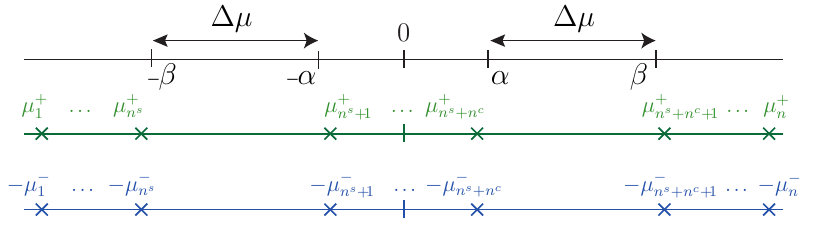}
\caption{Spectra of forward and backward FTLEs illustrating the gaps.} \label{slow_spectra_fig}
\end{figure}

 In Def.~\ref{definition1}, Property 1
 ensures that common gaps in
the forward and backward Lyapunov spectra not only exist, but also
separate the spectra in a dimensionally consistent manner, a relaxed version of Lyapunov regularity \cite{Pesin}. The
consistency between the forward and backward spectra is illustrated in Fig.~\ref{slow_spectra_fig}
 where the bounds and forward and backward exponents
 are plotted on aligned different copies of the real line for clarity.
 The exponents for particular values
 of $T$ and $\mathbf{x}$ are pictured, but note that Property 1 requires this structure for all $(T,
\mathbf{x})\in \mathcal{T}\times \mathcal{X}$. The symmetry of the gaps with respect to zero is not necessary but is
assumed here to simplify the presentation.
 The use of times up to $\overline{T}$ means that the computation of the Lyapunov exponents and vectors
 involves trajectories which, though they begin in $\mathcal{X}$,
 extend (unless $\mathcal{X}$ is $\phi-$invariant) into the larger set
 \begin{equation}\label{Kbar}
 \begin{aligned}
 \mathcal{X}_{ext}:=\{\mathbf{y}\in\mathbb{R}^n :\hspace{0.5mm} &\mathbf{y}=\phi(t,\mathbf{x})
		      \quad\textrm{or} \quad\mathbf{y}=\phi(-t,\mathbf{x})\\
		      &\quad\textrm{for some} \quad (t,\mathbf{x})\in \mathcal{T}\times \mathcal{X} \}.
 \end{aligned}
 \end{equation}
 The delayed start time $t_s$ provides a grace period
over which the FTLE bounds do not have to be satisfied,
in order to accommodate ``non-modal" behavior \cite{Schmid}; see Section 7.1 for a clarifying example.  $\overline{T}$
is the largest common time over
which the uniformity in the exponents holds. We note that because $\overline{T}$ must apply at each $\mathbf{x}$, for a
particular $\mathbf{x}$ larger forward and backward averaging times  may be possible; this property is exploited for the
example in Section \ref{MCKexample}. Given viable $\Delta\mu$, $t_s$, and $\overline{T}$, it can be
stated that the Lyapunov subspaces are resolvable over at least $\Delta\mu (\overline{T}-t_s)$
convergence time constants.

In Property 2, the subspaces $\mathcal{E}^s(\overline{T},\mathbf{x})$, $\mathcal{E}^c(\overline{T},\mathbf{x})$
and $\mathcal{E}^u(\overline{T},\mathbf{x})$ must uniformly define a splitting of the tangent space -- a finite-time
version of Oseledec's decomposition \cite{Pesin,Oseledec}. This condition is a transversality
requirement. The continuity of the splitting follows from the continuous dependence of $\Phi(\overline{T},\mathbf{x})$
on $\mathbf{x}$. We focus on the subspaces for $\overline{T}$ for the following reason. If the hypotheses of
Proposition \ref{theorem1} were applicable and the $T\rightarrow
\infty$
limits could be computed,
  then we could compute the forward and backward Lyapunov subspaces at each point of ${\mathcal{X}}$ for
arbitrarily large averaging times $T$ and  these subspaces would
converge to form $\Phi$-invariant subbundles \cite{Pesin}. Limited to $T\leq\overline{T}$ we should
use $T=\overline{T}$ to obtain subspaces that approximate the \textit{ideal} invariant subspaces
as closely as possible within the available averaging times. An argument similar to that in the proof of
Proposition \ref{theorem1} can be used to show that ${\cal
L}^+_{n^s}(T,\mathbf{x})$, ${\cal L}^-_{n^s+1}(T,\mathbf{x})$,
${\cal L}^+_{n^s+n^c}(T,\mathbf{x})$,
 and ${\cal L}^-_{n^s+n^c+1}(T,\mathbf{x})$ approach, with increasing $T$, fixed subspaces at least at a rate
proportional to $e^{-{\Delta\mu} T}$, and consequently so do the subspaces $\mathcal{E}^s(\overline{T},\mathbf{x})$,
$\mathcal{E}^c(\overline{T},\mathbf{x})$
and $\mathcal{E}^u(\overline{T},\mathbf{x})$.

Although the bounds $\alpha$ and $\beta$ give some indication of the timescales, in Property 3, the action of the transition matrix on vectors in the particular subspaces of the splitting in Property 2 is
characterized by exponential bounds. A procedure for determining $\nu$ and $\sigma$ is given in Section
\ref{diagnosis}. The time interval $\mathcal{T}_c$ over which the bounds apply is truncated at both ends. The delayed start time avoids the non-modal behavior and the cut-off time $t_c$ avoids a potential final transient from $t_c$ to $\overline{T}$ where a subspace
rotates away
from
the ideal asymptotic subspace it is intended to approximate. For a two-timescale set, $\nu -
\sigma$ is only
required to be
positive, but see the interpretation in Subsection \ref{Sig}.

\subsection{Invariant Manifold Approximation}

If ${\cal X}$, now assumed to be a domain of $\mathbb{R}^n$, is a finite-time uniform two-timescale set, we postulate a
corresponding manifold structure for the flow of the nonlinear system (\ref{nldyn}). The characteristics of the
two-timescale set provide the dimensions and orientations of the manifolds. In particular, in the general case where
none of the dimensions $n^s$, $n^c$, or $n^u$ is zero, center $\mathcal{W}^c$, center-stable $\mathcal{W}^{cs}$, and
center-unstable $\mathcal{W}^{cu}$ invariant manifolds can be postulated along with a corresponding invariant splitting
$T_x\mathcal{X} = \mathcal{E}^s(\mathbf{x})\oplus \mathcal{E}^c(\mathbf{x}) \oplus \mathcal{E}^u(\mathbf{x})$. Points on
the postulated invariant manifolds satisfy the conditions
\begin{equation}%\label{finiteS}
\begin{array}{ccl}
\mathbf{x}\in\mathcal{W}^c \Rightarrow \mathbf{f}(\mathbf{x})\in \mathcal{E}^c(\mathbf{x})\\
\mathbf{x}\in\mathcal{W}^{cs} \Rightarrow \mathbf{f}(\mathbf{x})\in \mathcal{E}^s(\mathbf{x})\oplus
\mathcal{E}^c(\mathbf{x})\\
\mathbf{x}\in\mathcal{W}^{cu} \Rightarrow \mathbf{f}(\mathbf{x})\in \mathcal{E}^c(\mathbf{x}) \oplus
\mathcal{E}^u(\mathbf{x})
\end{array}\end{equation}

\noindent Approximating the postulated invariant splitting with our finite-time non-invariant splitting, we  can search for points that satisfy the subspace membership conditions (which will be posed as orthogonality conditions in the
next section). This leads to the definition of a finite-time center manifold.
 \begin{definition} Given a uniform finite-time two-timescale set $\mathcal{X}$, a
 finite-time center manifold is an $n^c$-dimensional submanifold of $\mathcal{X}$
 denoted $\mathcal{W}^c(\overline{T})$ such that $\mathbf{f}(\mathbf{x})\in \mathcal{E}^c(\overline{T},\mathbf{x})$
  for all $\mathbf{x}\in \mathcal{W}^c(\overline{T})$.\end{definition}

\noindent Analogous definitions can be given for the finite-time manifolds $\mathcal{W}^{cs}(\overline{T})$ and
$\mathcal{W}^{cu}(\overline{T})$.

\subsection{Interpretation and Significance}\label{Sig}
Consider the scenario in which the behavior of a system $\dot{\mathbf{x}}=\mathbf{f}(\mathbf{x})$
over the time interval  $[0, t_f]$ is of interest. Assume set $\mathcal{X}$ covers the region of state space in which this behavior takes place and has been diagnosed a uniform finite-time
two-timescale set with time constants $\nu^{-1}$ and $\sigma^{-1}$. If $t_f$ is
much larger than $\nu^{-1}$ and smaller than or similar to $\sigma^{-1}$, then there is slow-fast behavior in the
tangent linear dynamics relative to the time interval of
interest.\footnote{If there is more than one way to separate the FTLE spectra to satisfy Def. \ref{definition1}, then
the value of $t_f$ of interest can suggest which way to consider.} Further, if $t_s$ and $\overline{T}-t_c$ are  small
fractions of $t_f$, then the exponential bounds apply
on most of the time interval. If an $n^c$-dimensional invariant center manifold is present in $\mathcal{X}$ and $n^u=0$ ($n^s=0$), trajectories in a neighborhood of the manifold will approach, during a small
fraction of $t_f$, the center
manifold in forward (backward) time, and one could approximate the behavior over most of the time interval as
behavior on the reduced-order manifold. If both $n^u$ and $n^s$ are nonzero, then trajectories will approach the invariant
center-unstable
$\mathcal{W}^{cu}$ (the invariant center-stable $\mathcal{W}^{cs}$ ) manifold in forward (backward) time; the example in Section \ref{3DnonlinEx} illustrates this.

Our interest in a center manifold in $\mathcal{X}$ can be related to interest in a fixed point, in the sense that locating each and analyzing the linearized dynamics provides valuable information about the flow. At a fixed point there is complete equilibrium, whereas on a center manifold there is partial equilibrium. FTLA characterizes the behavior of trajectories in a neighborhood of a normally hyperbolic center manifold, analogous to the role of linear analysis at a hyperbolic fixed point. The non-modal behavior that, for the purpose of locating the center manifold, has been excluded from influencing the FTLEs, via the delayed start time $t_s$, can affect the size of the neighborhood \cite{Schmid} and should be considered.

In the remainder of the paper, we focus on computing points on finite-time center manifolds.  For the general normally hyperbolic case, this
requires obtaining $[\mathcal{E}^c]^\bot$ by intersecting filtrations from forward and backward integration of the tangent linear dynamics. Points on $\mathcal{W}^{cu}$ and
$\mathcal{W}^{cs}$ can also be determined using subspace membership conditions and can benefit the solution of certain
boundary-value problems \cite{ACC09, Guck09, Rao99, Rao00, Kopell}. As mentioned earlier, center manifolds need not be slow manifolds. At a point $\mathbf{x}\in W^c(\overline{T})$, the exponential bounds for $\mathcal{E}^c(\overline{T},\mathbf{x})$ constrain the rate of
change in the length of $\mathbf{f}(\mathbf{x})$ but the FTLA characterization does not constrain the length of $\mathbf{f}(\mathbf{x})$ and leaves rotational freedom, even fast rotation is possible.

\section{Finite-Time Two-Timescale Set and Center Manifold - Procedure} \label{procedure}
 If
the goal is only to diagnose two-timescale behavior and determine
the tangent space structure, then $\mathcal{X}$ can be any subset of
$\mathbb{R}^n$. If one also
wants to search for a center manifold, then  $\mathcal{X}$ is typically
a domain, or a set of grid points on a domain, because it will be necessary to iteratively search for
points that satisfy center manifold conditions in a state space region
of full dimension. As mentioned, simulation experience with a set of boundary conditions of interest
could suggest the domain $\mathcal{X}$ to explore. Establishing \textit{a priori} that one is searching for a center manifold in a two-timescale set ensures that a uniform splitting exists and can be resolved; however it is also possible to proceed directly to the search and assess the uniformity of the timescales and splitting in the process.

\subsection{Diagnosing a Finite-Time Two-Timescale
Set}\label{diagnosis} The three properties in Def. \ref{definition1} are checked on
$\mathcal{X}$. To check Property 1, FTLEs are
computed  to determine if there exist a $t_s$ and $\overline{T}$ for which there is a
pattern as illustrated in Fig.~\ref{slow_spectra_fig}
uniformly in $\mathbf{x}$ and for all $T\in \mathcal{T}$. Regarding
uniformity, the individual exponents can vary with $T$ and
$\mathbf{x}$ as long as there is a sufficiently large uniform gap.
 However, unless $\mathcal{X}$ is
$\phi-$invariant, the set $\mathcal{X}_{ext}$ (see (\ref{Kbar}))
grows with $T$; an upper limit on $T$ may be required to avoid non-uniform behavior. If a uniform spectral gap is
present, then the appropriate values of the constants $n^s$, $n^c$, $n^u$, $t_s$, $t_c$, $\overline{T}$ and $\Delta\mu$
are determined. Based on the initial survey, $\mathcal{X}$ could be adjusted.

If the tangent space structure is resolvable over at least 3-5 convergence time constants, then the subspaces
$\mathcal{E}^s(\overline{T},\mathbf{x})$, $\mathcal{E}^c(\overline{T},\mathbf{x})$ and
$\mathcal{E}^u(\overline{T},\mathbf{x})$ are constructed and Property 2 is checked. The dimensions of these subspaces
sum to $n$, but each pair of subspaces must intersect transversely to provide the splitting. We note that the
convergence of the subspaces can be checked
directly by monitoring the distance between the subspaces with
increasing averaging time (illustrated in Section \ref{examples}).

A means \cite{PesinHB} of confirming that an invariant splitting exists close to the splitting
(\ref{FTsplitting}) uses cones defined as follows.  The {\em cone} at $\mathbf{x}\in \mathbb{R}^n$ centered on
the subspace $S(\mathbf{x})\subset T_\mathbf{x}\mathbb{R}^n$ with angle $\psi\in(0,\pi/2)$ is given by
\begin{equation}\label{ConeDef} C(\mathbf{x},S(\mathbf{x}),\psi):=\{\mathbf{v}\in T_\mathbf{x}\mathbb{R}^n:
\angle(\mathbf{v},S(\mathbf{x}))<\psi\},
\end{equation}
where $\angle(\mathbf{v},S(\mathbf{x}))$ is the angle between $\mathbf{v}$ and its orthogonal projection in
$S(\mathbf{x})$. One tries  to verify that there are families of stable, unstable, center-stable and
center-unstable cones
\begin{equation}\nonumber
\begin{aligned}
&C^s(\mathbf{x},\psi)=C(\mathbf{x},\mathcal{E}^s(\overline{T},\mathbf{x}),\psi), \hspace{5mm}
 C^u(\mathbf{x},\psi)=C(\mathbf{x},\mathcal{E}^u(\overline{T},\mathbf{x}),\psi),\\
&C^{cs}(\mathbf{x},\psi)=C(\mathbf{x},\mathcal{E}^{cs}(\overline{T},\mathbf{x}),\psi), \hspace{3mm}
 C^{cu}(\mathbf{x},\psi)=C(\mathbf{x},\mathcal{E}^{cu}(\overline{T},\mathbf{x}),\psi),
\end{aligned}
\end{equation} where
\begin{equation}\nonumber
\mathcal{E}^{cs}(\overline{T},\mathbf{x})=\mathcal{E}^c(\overline{T},\mathbf{x})\oplus
\mathcal{E}^s(\overline{T},\mathbf{x}), \hspace{2.5mm}
\mathcal{E}^{cu}(\overline{T},\mathbf{x})=\mathcal{E}^c(\overline{T},\mathbf{x})\oplus
\mathcal{E}^u(\overline{T},\mathbf{x}),
\end{equation}
such that
\begin{equation}\nonumber
\setlength{\abovedisplayskip}{0pt}
\begin{aligned}
&\Phi(-t,\mathbf{x}) C^s(\mathbf{x},\psi)\subset C^s(\phi(-t,\mathbf{x}),\psi), \\
&\Phi(t,\mathbf{x}) C^u(\mathbf{x},\psi)\subset C^u(\phi(t,\mathbf{x}),\psi),\\
&\Phi(-t,\mathbf{x}) C^{cs}(\mathbf{x},\psi)\subset C^{cs}(\phi(-t,\mathbf{x}),\psi),\\
&\Phi(t,\mathbf{x}) C^{cu}(\mathbf{x},\psi)\subset C^{cu}(\phi(t,\mathbf{x}),\psi),
\end{aligned}
\end{equation}
 for all $t\in \mathcal{T}$. The notation $\Phi(-t,\mathbf{x}) C^s(\mathbf{x},\psi)$ means the subspace at
$\phi(-t,\mathbf{x})$ are obtained by backward propagation of all the vectors in $C^s(\mathbf{x},\psi)$. For the cone
conditions to be satisfied, $\psi$ must be large enough that the cones
contain the actual invariant subspaces. The size of $\psi$ could be  iteratively reduced  to get an estimate of how
close the splitting is to being invariant.

For each $\mathbf{x}\in\mathcal{X}$, the subspaces (\ref{EdefsT}) that define the splitting of the tangent space
$T_\mathbf{x}{\mathbb{R}^n}$ at $\overline{T}$ can be expressed as the column spans (i.e., range spaces) of the
following matrices
\begin{equation}\label{subsp_Tbar}
 \begin{aligned}
  E^s(\overline{T},\mathbf{x}) &=
[\mathbf{l}^+_1(\overline{T},\mathbf{x}),...,\mathbf{l}^+_{n^s}(\overline{T},\mathbf{x})],\\
  E^c(\overline{T},\mathbf{x}) &= null\left[([E^c(\overline{T},\mathbf{x})]^{\perp})^T\right]
,\\
  E^u(\overline{T},\mathbf{x}) &=
[\mathbf{l}^-_{n^s+n^c+1}(\overline{T},\mathbf{x}),...,\mathbf{l}^-_{n}(\overline{T},\mathbf{x})],\\
 \end{aligned}
\end{equation}
\noindent
$[E^c(\overline{T},\mathbf{x})]^{\perp}$ is given in terms of the FTLVs in the next subsection. We have used `$null(M)$' to denote the mapping from matrix $M$ to
an orthonormal matrix whose column span is the null space of the
matrix $M$.

To check Property 3, we check if $\nu > \sigma$ after computing the constants $\nu$ and $\sigma$ as
\begin{equation}\nonumber
 \nu=\min\{-\overline{\mu}^s,\underline{\mu}^u,\underline{\mu}^s,-\overline{\mu}^u \} \hspace{1mm},\hspace{5mm}
 \sigma=\max\{|\underline{\mu}^{c+}|,|\overline{\mu}^{c+}|,|\overline{\mu}^{c-}|,|\underline{\mu}^{c-}|\}\\
\end{equation}
where
\begin{equation}\label{mu_lim}
\begin{aligned}
 \overline{\mu}^s=\sup\limits_{(T,\mathbf{x})\in
\mathcal{T}_c\times\mathcal{X}}\mu^{s+}_{n^s}(T,\mathbf{x}),
\hspace{5mm} &
\underline{\mu}^{c+}=\inf\limits_{(T,\mathbf{x})\in
\mathcal{T}_c\times\mathcal{X}}\mu^{c+}_{1}(T,\mathbf{x}),\\
\underline{\mu}^u=\inf\limits_{(T,\mathbf{x})\in
\mathcal{T}_c\times\mathcal{X}}\mu^{u+}_{1}(T,\mathbf{x}),
\hspace{5mm} &
\overline{\mu}^{c+}=\sup\limits_{(T,\mathbf{x})\in
\mathcal{T}_c\times\mathcal{X}}\mu^{c+}_{n^c}(T,\mathbf{x}),\\
\underline{\mu}^s=\inf\limits_{(T,\mathbf{x})\in
\mathcal{T}_c\times\mathcal{X}}\mu^{s-}_{n^s}(T,\mathbf{x}),
\hspace{5mm} &
\overline{\mu}^{c-}=\sup\limits_{(T,\mathbf{x})\in
\mathcal{T}_c\times\mathcal{X}}\mu^{c-}_{1}(T,\mathbf{x}),\\
\overline{\mu}^u=\sup\limits_{(T,\mathbf{x})\in
\mathcal{T}_c\times\mathcal{X}}\mu^{u-}_{1}(T,\mathbf{x}),
\hspace{5mm} &
\underline{\mu}^{c-}=\inf\limits_{(T,\mathbf{x})\in
\mathcal{T}_c\times\mathcal{X}}\mu^{c-}_{n^c}(T,\mathbf{x}).\\
\end{aligned}
\end{equation}

\noindent The FTLEs for each subspace as needed in (\ref{mu_lim}) are computed as
\begin{equation}\nonumber
\mu^{j\pm}_i(T)=\frac{1}{T}\ln\left(\Sigma^{j\pm}_{ii}\right)
\hspace{5mm} i=1,...,n^j, j=s,c,u,\\
\end{equation}
where the diagonal matrices $\Sigma^{j\pm}$ are obtained from the
singular value decompositions
\begin{equation}\label{SVD_prop}
N^{j\pm}(\pm T,\phi(\pm T,\mathbf{x}))\cdot \Sigma^{j\pm}(\pm T, \mathbf{x})\cdot L^{j\pm}(\pm T,\mathbf{x})
=\Phi(\pm T,\mathbf{x}) E^{j}(\overline{T},\mathbf{x})\\
\end{equation}
\noindent and the subscript `$ii$' on $\Sigma$ denotes the $i^{th}$ diagonal element of that matrix.

\subsection{Computing Points on a Finite-Time Center Manifold} Provided that $\mathcal{X}$ satisfies
 Def.~\ref{definition1}, where ${\cal X}$ is now assumed to be a domain of
$\mathbb{R}^n$, we can
 look for a normally hyperbolic center manifold in ${\cal
X}$. Within ${\cal X}$, the points in the set
\begin{equation}\label{finiteS}
\{\mathbf{x}\in\mathcal{X}:\langle\mathbf{f}(\mathbf{x}),\mathbf{w}\rangle=0,
\forall \mathbf{w}\in [\mathcal{E}^c(\overline{T},\mathbf{x})]^\perp\}\end{equation}
satisfy a necessary condition for being on a
finite-time center manifold. Whether or not
this set, or a subset of it, is a submanifold of ${\cal X}$ has to be determined to the extent it can from numerical results, which is
also the case with the ILDM method \cite{MaasPope}. We proceed under the assumption that
the set is a manifold that can locally be parametrized by $n^c$ of the $n$ system coordinates and
represented as a graph. Situations with folded or multiple center manifolds, or where the set is an object of
fractal dimension are possible; see \cite{Lorenz} for an example of analyzing the case of a fractal attractor, though
not characterized as the set (\ref{finiteS}).

Rather than use
eigenvectors of $D\mathbf{f}(\mathbf{x})$ to form an approximate
basis for the orthogonal complement to $\mathcal{E}^c$ as in the ILDM method \cite{MaasPope},
we use the appropriate Lyapunov vectors to form the approximate basis as prescribed in the following proposition.

\begin{proposition}\label{proposition3}
 On a uniform finite-time two-timescale set $\mathcal{X}$, at each $\mathbf{x}$, the vectors
\begin{equation}\label{prop3eq}
\begin{aligned}
\mathbf{l}^-_1(\overline{T},\mathbf{x}),\dots,\mathbf{l}^-_{n^s}(\overline{T},\mathbf{x}),
\mathbf{l}^+_{n^s+n^c+1}(\overline{T},\mathbf{x}),\dots,\mathbf{l}^+_{n}(\overline{T},\mathbf{x})
\end{aligned}
\end{equation}
form a basis for $[\mathcal{E}^c(\overline{T},\mathbf{x})]^\bot$.
\end{proposition}
{\it Proof:} In Def. \ref{definition1}, Property 2, the $n^c$-dimensional center subspace is given by $\mathcal{E}^c(\overline{T},\mathbf{x})=$
${\cal L}^+_{n^s+n^c}(\overline{T},\mathbf{x})\cap {\cal L}^-_{n^s+1}(\overline{T},\mathbf{x})$. Using an identity from \cite{Isidori}, we
have $[\mathcal{E}^c(\overline{T},\mathbf{x})]^\bot=[{\cal L}^-_{n^s+1}(\overline{T},\mathbf{x})]^\bot ~\oplus$ $
[{\cal L}^+_{n^s+n^c}(\overline{T},\mathbf{x})]^\bot$. The proposition then follows from the facts:
$[{\cal L}^-_{n^s+1}(\overline{T},\mathbf{x})]^\bot= span~\{\mathbf{l}^-_1(\overline{T},\mathbf{x}),\dots,\mathbf{l}^-_{n^s}(\overline{T},\mathbf{x})\}$
and $[{\cal L}^+_{n^s+n^c}(\overline{T},\mathbf{x})]^\bot=span~\{\mathbf{l}^+_{n^s+n^c+1}(\overline{T},\mathbf{x}),\dots,
\mathbf{l}^+_{n}(\overline{T},\mathbf{x})\} .$  $\hfill\blacksquare$

The set (\ref{finiteS}) is thus the solution set for the system of orthogonality conditions
\begin{equation}\label{OrthogConds}
 \begin{aligned}
  \langle \mathbf{f}(\mathbf{x}),\mathbf{l}^-_i(\overline{T},\mathbf{x})\rangle &= 0, i=1,\dots,n^s\\
  \langle \mathbf{f}(\mathbf{x}),\mathbf{l}^+_j(\overline{T},\mathbf{x})\rangle &= 0, j=n^s+n^c+1,\dots,n\\
 \end{aligned}
\end{equation}
\noindent In order to obtain solutions of the algebraic equations, we
designate $n^c$ components of $\mathbf{x}$ as independent variables to parameterize the manifold and determine the values of the remaining $n-n^c$ components, the dependent variables, that satisfy the orthogonality conditions in (\ref{finiteS}). Because the FTLVs are in numerical form, we use a successive approximation approach described in Section \ref{MCKexample}.
 This is repeated for a grid on the space of independent variables. The directions of the Lyapunov vectors indicate how
to
separate the coordinates of $\mathbf{x}$ into independent and
dependent variables, i.e., how to locally parametrize the postulated ${\cal W}^c(\overline{T})$. The independent variables must be
chosen
such that their coordinate axes are not parallel to any directions in $[\mathcal{E}^c]^\bot$. Different independent
variables might be required for different sections of the center manifold.

Consideration of the planar system
 \begin{equation}
 \dot{\mathbf{x}}=\mathbf{f}(\mathbf{x})=\mathbf{e}^c(x_1)g(x_1)+\mathbf{e}^s(x_1)h(x_1,x_2)
 \end{equation}
 where $\mathbf{e}^c$ and $\mathbf{e}^s$ are unit basis vectors for the exact center and stable subspaces and $g$ and $h$ are scalar functions, provides some insight into what the accuracy of the finite-time manifold depends on. Assume the there is an invariant center manifold given by the solution set for $h(x_1,x_2)=0$ and parametrizable by $x_1$. For a particular value of $x_1$, let unit vector $\mathbf{w}$ be the approximation of the direction orthogonal to $\mathbf{e}^c$ at the correct location $x_2^c$ of the manifold. Assume $<\mathbf{w},\mathbf{e}^c>=\sin(\varepsilon)$ and $<\mathbf{w},\mathbf{e}^s>=\sin(\delta+\varepsilon)$, where $\varepsilon$ is the angular error in approximating $\mathbf{e}^c$ and $\delta$ is the angle between $\mathbf{e}^c$ and $\mathbf{e}^s$. Then the manifold error is
 \begin{equation}
 \Delta x_2=-\frac{\sin(\varepsilon)}{\sin(\delta+\varepsilon)}g(x_1)\left(\frac{\partial{h}}{\partial{x_2}}(x_1)\right)^{-1}
 \end{equation}
 In addition to the subspace error $\varepsilon$, the manifold error depends on the angle $\delta$ between $\mathbf{e}^c$ and $\mathbf{e}^s$ and the ratio $g/(\partial{h}/\partial{x_2})$.

\subsection{Numerical Methods for FTLA} \label{nummeths} Numerical methods for FTLA are addressed in
\cite{Adrover,adrover06,Dieci,Geist,Samelson} and the references therein. For completeness, the
 methods used for the computations presented in the next
section are described in this subsection. All the computations are
done in the Matlab$^{\circledR}$ environment. The numerical
integration of the nonlinear state equations and the corresponding
linear variational equations is performed with the `ode45'
integrator.

The FTLEs and FTLVs associated with an initial state $\mathbf{x}$
are computed for an averaging time $T$ either by SVD or QR
factorization. Only the computation of the forward-time FTLE/Vs is
described, since the computation of the backward-time FTLE/Vs is
analogous. The first step of both methods is to integrate the
nonlinear state equations from $t=0$ to $t=T$ and save the values of
$\phi(t,\mathbf{x})$ at the $N$ equally spaced times $\Delta t,
2\Delta t, \dots, N\Delta t$, where $N\Delta t=T$.

In the SVD method, the transition matrix is computed and then the
SVD is applied. The transition matrix is computed by integrating,
simultaneously, the nonlinear equations and the associated linear variational
equations over each segment of the base space trajectory, with the
state initialized with the saved value at the beginning of the
segment and the transition matrix initialized with the identity
matrix. Using the notation $\Phi_k^{\Delta t}=\Phi(\Delta
t,\phi[(k-1)\cdot \Delta t, \mathbf{x}])$ for $k=1,2,...,N$, the
transition matrix is constructed from the transition matrices for
the segments as $\Phi(T,x)=\Phi_N^{\Delta t}\cdots\Phi_2^{\Delta
t}\Phi_1^{\Delta t}$. The resulting transition matrix is then
factored as $\Phi(T,\mathbf{x})=N^+\Sigma^+ (L^+)^T$ using the `svd'
command in Matlab$^{\circledR}$. Each FTLE
is obtained by $\mu^+_i(T,\mathbf{x})=\frac{1}{T}\ln{\sigma^+_i}$,
where $\sigma_i$ is the $i^{th}$ singular value of $\Phi$, the
positive square root of the $i^{th}$ diagonal element of $\Sigma^+$.
If this procedure does not produce FTLEs in the ascending order we
have assumed in our notation, the FTLEs and associated FTLVs are rearranged to conform. The FTLVs $\mathbf{l}_i^+(T,
\mathbf{x}), i=
1,\dots, n$ are the column vectors of $L^+$.

For a given trajectory from $\mathbf{x}$ to $\phi(T, \mathbf{x})$,
for a particular $T$, we have the option of computing the
Lyapunov vectors at $\mathbf{x}$ and at $\phi(T, \mathbf{x})$ by
forward or backward integration. Because
$\Phi(-T,\phi(T,\mathbf{x}))=\Phi^{-1}(T,\mathbf{x})$, it follows
that $L^+(T,\mathbf{x})=N^-(T,\phi(T,\mathbf{x}))$ and
$N^+(T,\mathbf{x})=L^-(T,\phi(T,\mathbf{x}))$. As pointed out by others, e.g. in \cite{Legras}, it is best to compute
$L^+(T,\mathbf{x})$ by backward integration from $\phi(T, \mathbf{x})$ and $L^-(T,\phi(T,\mathbf{x}))$ by forward
integration from $\mathbf{x}$ so that the vectors and subspaces  one is seeking are those to which the linear flow
naturally carries the vectors and subspaces. The QR method is based on this strategy.

In the QR method, a segmented approach is also used
\cite{Dieci}. For the $k^{th}$ segment, after the transition matrix
is computed as described in the previous paragraph, the $Q_{k-1}$
matrix associated with the state at the end of the previous segment
is propagated by the transition matrix to the end of the
$k^{th}$ segment and the $Q_kR_k$ factorization of the resulting
matrix is obtained, as summarized by
\begin{equation}\label{eqn:QRstep}
\Phi_k^{\Delta t}Q_{k-1}=Q_kR_k.
\end{equation}
This sequence of operations for $k=1,\dots,N$ must be initialized by
prescribing $Q_o$; typically the identity matrix is used
\cite{Dieci,Geist}. It then follows that
\begin{equation}\label{eqn:QR}
\Phi(T,\mathbf{x})Q_o=Q(T,\mathbf{x})R
\end{equation}
where $Q(T,\mathbf{x})=Q_N$ and $R=R_NR_{N-1}\dots R_2R_1$. For
almost every $Q_o$, as $T$ increases, $Q(T,\phi(T,\mathbf{x}))$ will approach
$N^+(T,\phi(T,\mathbf{x}))$ and the diagonal elements of $R$ will approach the diagonal elements of $\Sigma^+$ in the
absence of numerical errors. Note that, for any $T$,  if we choose $Q_o=L^+(T,\mathbf{x})$, then
$Q(T,\mathbf{x})=N^+(T,\mathbf{x})$, or equivalently
$Q(T,\mathbf{x})=L^-(T,\phi(T,\mathbf{x}))$, and $R=\Sigma^+$. In our experience, the QR method is generally more
reliable than the SVD method for calculating the FTLE/Vs for longer averaging times. For shorter times, as needed to
compute the exponent bounds, the SVD should be used.

\section{Application Examples} \label{examples} Several application examples are
presented to demonstrate the use of the FTLA methodology. The first example provides insight into the start and cut-off
times used in Definition \ref{definition1}. The angles between the relevant vectors and subspaces are intentionally small to illustrate how the FTLA method handles the consequences. The remaining three examples illustrate the FTLA methodology for 2D, 3D, and
4D systems, the second two involving normally hyperbolic center manifolds. For these examples, the subspaces in the splitting are are separated by angles of at least 45 degrees, and the center manifolds are slow manifolds. Given that our initial motivation for
developing the FTLA method for determining a slow manifold was to improve the accuracy of the ILDM method in situations
where the ILDM method is known to be inaccurate \cite{Kaper}, the FTLA method results are compared to the results
obtained with the ILDM method. The comparisons thus focus on the differences between using FTLA and analyzing the
tangent linear dynamics as if they were time-invariant and how these differences map into manifold errors. Note that other means of improving the results of the ILDM
method have been developed, e.g., \cite{NafeMaas}.

\subsection{Example for Understanding Start and Cut-Off Times }

Properties 1 and 3 in Def. \ref{definition1} involve truncating the time
interval at the beginning and end,
using the start time $t_s$ and the cut-off time $t_c$. The initial
transient behavior that is excluded is
associated with coordinate-dependent angles between certain vectors
within the ideal asymptotic stable, center and unstable subspaces toward which the finite-time subspaces are converging. The
final transient behavior that is excluded is produced by the lack of
$\Phi$-invariance of the finite-time subbundles $\mathcal{E}^s$,
$\mathcal{E}^c$ and $\mathcal{E}^u$. To illustrate the behaviors and the roles of
the constants $t_s$ and $t_c$, we
consider a 7D system, $\mathbf{\dot{x}}=\mathbf{f}(\mathbf{x})$, at an
equilibrium point $\mathbf{x}_e$, i.e., for
$\mathcal{X}=\{\mathbf{x}_e\}$, with
\begin{equation}D\mathbf{f}(\mathbf{x}_e)=
\begin{bmatrix} -5.4 &   1  &   0  &   0  &   0 &   0  &   0\\
          0  & -5.2 &   0  &   0  &   30 &   0 &   0\\
          0  &   0  & -0.3 &   0  &   0  &   0 &   10\\
          0  &   0  &   0  & -0.1 &   0  &   0 &   0\\
          0  &   0  &   0  &   0  &  0.2 &   0 &   0\\
          0  &   0  &   0  &   0  &   0  &  4.0 &   8\\
          0  &   0  &   0  &   0  &   0  &   0 &  4.6\\
\end{bmatrix}.
\end{equation}
The triangular form of $D\mathbf{f}(\mathbf{x}_e)$ allows simple control
of the timescales, the important angles, and
the degree of dynamic
coupling via specification of the diagonal and off-diagonal elements.

Barring numerical errors, in the limit $\overline{T}\rightarrow\infty$,
the FTLEs will converge to the eigenvalues of
$D\mathbf{f}(\mathbf{x}_e)$, i.e., the diagonal elements, and the
subspaces $\mathcal{E}^s$, $\mathcal{E}^c$ and
$\mathcal{E}^u$ will converge to the stable, center and
unstable eigenspaces, i.e., the subspaces spanned by the appropriate
subset of the eigenvectors of
$D\mathbf{f}(\mathbf{x}_e)$ -- the stable
eigenspace spanned by the eigenvectors for the eigenvalues
$\lambda_1=-5.4$ and $\lambda_2=-5.2$, the center eigenspace spanned by
the eigenvectors for
$\lambda_3=-0.3$, $\lambda_4=-0.1$ and $\lambda_5=0.2$ and  the unstable
eigenspace spanned by the eigenvectors for
$\lambda_6=4.0$ and $\lambda_7=4.6$.

In order to determine the cut-off time $t_c$, the FTLEs for the subspaces
$\mathcal{E}^s$, $\mathcal{E}^c$ and $\mathcal{E}^u$ computed for a
finite $\overline{T}$ ($\mu^{j\pm}_i, i=1,...,n^j,
j=s,c,u$) are considered to determine the exponential bounds as described in Section \ref{diagnosis}. For sufficiently large
finite $\overline{T}$, the subspaces
$\mathcal{E}^s$, $\mathcal{E}^c$ and $\mathcal{E}^u$ will closely
approximate the corresponding fixed eigenspaces,
but when propagated to $\overline{T}$, there will be a final
boundary-layer in which the subspaces rotate away from the
eigenspaces and this will affect the behavior of the FTLEs. For example,
the stable eigenspace is asymptotically stable
in backward time and unstable in forward time with respect to
neighboring equi-dimensional subspaces. Thus, when propagated forward in time, the
finite-time approximation $\mathcal{E}^s$ will rotate away from the stable eigenspace. This is a
non-uniform rotation taking place primarily near the
time $\overline{T}$ for which $\mathcal{E}^s$ was
computed. In general $\mathcal{E}^s$ and
$\mathcal{E}^c$ will rotate toward
$\mathcal{E}^u$ in forward time, and $\mathcal{E}^c$ and $\mathcal{E}^u$
will rotate toward $\mathcal{E}^s$ in backward
time. The FTLEs for $\mathcal{E}^s$, $\mathcal{E}^c$ and
$\mathcal{E}^u$ will be similar to those for their eigenspace
counterparts except in cases involving  propagation in the
unstable
direction when the averaging time is near $\overline{T}$. Thus we
exclude a final transient
period long enough to avoid the corresponding deviations in the FTLEs.
Figure \ref{FTLE_EB_7DLTI} shows the backward and
forward FTLEs for each of the
three subspaces for $\overline{T}=6.0$. The final transients are short
and the deviations are
not large; the final transients that dictate $t_c=5.5$ are the ones for
the forward and backward propagations of the
center subspace. On the other hand, the determination of the start time $t_s$ comes from
the requirement of satisfying properties 1 and
3 of \ref{definition1}. Therefore, we will consider both the FTLEs that
define the
exponential bounds (\ref{EdefsT}) and the FTLEs that define the
spectral gap $\Delta\mu$. For convenience, we will refer
 to these two sets of FTLEs as $\mu^{\pm}_{EB}$ and $\mu^{\pm}$ respectively. The $\mu^{\pm}_{EB}$'s associated with a particular subspace, as functions
of $T$, will have an initial transient period,
if the subspace has dimension greater than one and there is one or more
pair of eigenvectors within the eigenspace
being approximated that are separated by an
angle less than $90^\circ$ in the coordinates being used. In
this example, the angles referred to are those between the eigenvectors
that span the stable, center, and unstable
eigenspaces. Angles less than $90^{\circ}$ are responsible for the
funnel-shaped initial transient behavior of the $\mu^{\pm}_{EB}$. For
instance, the angle between the two
eigenvectors associated with the two largest eigenvalues is $9.7^{\circ}$ and the backward $\mu^{-}_{EB}$
for $\mathcal{E}^u$ in the  $T\rightarrow 0$ limit (i.e. the opposites of the eigenvalues of the symmetric part of $D\mathbf{f}(\mathbf{x}_e)$,  which are
$-2.5$
and $-6.1$) are not consistent with the $\mu^{-}$ for most
averaging times up to $\overline{T}$; this is
referred to as non-modal behavior \cite{Schmid}. By excluding a period
$[0,t_s]$ the initial transient behavior is
eliminated. A similar argument can be made when considering the FTLEs
$\mu^{\pm}$.

%The following is a procedure for computing $t_s$. Once the constants $n^s$, $n^c$, $n^u$ and the
%cutoff time $t_c$ are determined, consider a point
%$\mathbf{x}\in\mathcal{X}$.
%(i)  Calculate $\Delta\mu_{EB}(\mathbf{x},t_c)$ and
%$\Delta\mu(\mathbf{x},\overline{T})$. These values are taken as
%references because at $t_c$ and $\overline{T}$ we are guaranteed to have
%a well-defined two-timescale separation and
%sufficient convergence of the subspaces $\mathcal{L}_i^{\pm}$.
%(ii) Multiply $\Delta\mu_{EB}(\mathbf{x},t_c)$ and
%$\Delta\mu(\mathbf{x},\overline{T})$ by a constant factor
%$0<\hat{k}<1$ and calculate the
%times ($t_s'(\mathbf{x})$, $t_s''(\mathbf{x})$) at which
%$\Delta\mu_{EB}(\mathbf{x},t)=\hat{k}\Delta\mu_{EB}(\mathbf{x},t_c)$ and
%$\Delta\mu_{EB}(\mathbf{x},t)=\hat{k}\Delta\mu(\mathbf{x},\overline{T})$.
%This guarantees uniformity of $\mu^{\pm}_{EB}$ and $\mu^{\pm}$ in
%$\mathcal{T}_c$ and $\mathcal{T}$ respectively. If one
%or both the previous equations do not have a solution, then set
%$t_s'(\mathbf{x})$ or $t_s''(\mathbf{x})$ or both to zero. Finally pick
%$t_s(\mathbf{x}) = \max(t_s'(\mathbf{x}), t_s''(\mathbf{x}))$.
%Repeat the procedure for all $\mathbf{x}\in\mathcal{X}$. The actual start time
%will then be
%$t_s=\sup_{\mathbf{x}\in\mathcal{X}} t_s(\mathbf{x})$.
%A larger value of $t_s$ allows tighter exponential bounds but shortens
%the time interval over which they apply.

\begin{figure}[ht]
\hspace{-2mm}
\includegraphics[scale=.36]{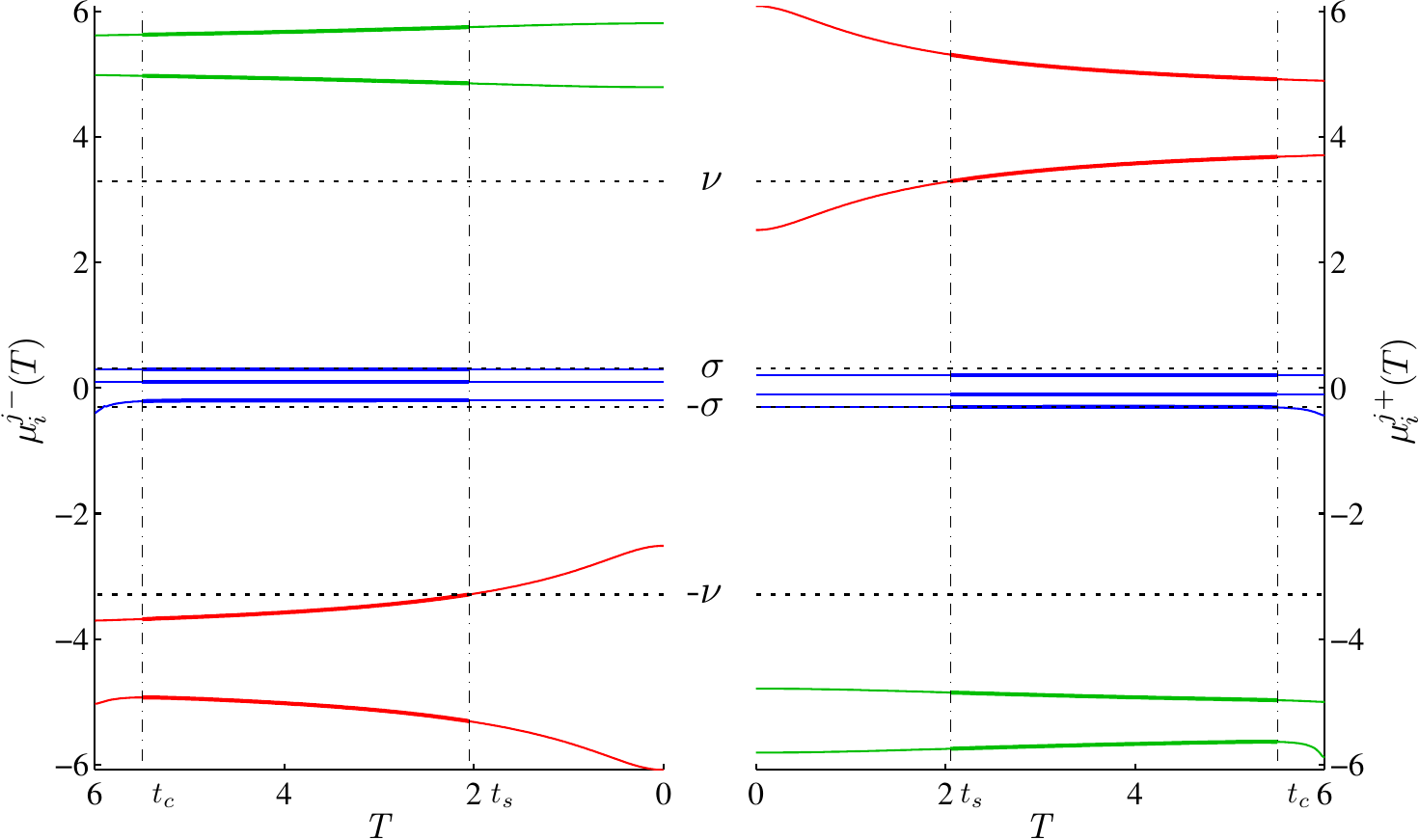}
\caption{Backward and forward FTLEs ($\mu^{j \pm}_{i}$ with $i=1,...,n^j$ and $j=s,c,u$) for the subspaces
$\mathcal{E}^s$ (green), $\mathcal{E}^u$ (red) and $\mathcal{E}^c$ (blue). The exponential bound constants $\sigma$ and
$\nu$ and the start and cutoff times $t_s$ and $t_c$ are shown.}
\label{FTLE_EB_7DLTI}
\end{figure}

\begin{figure}[ht]
\hspace{-2mm}
\includegraphics[scale=.355]{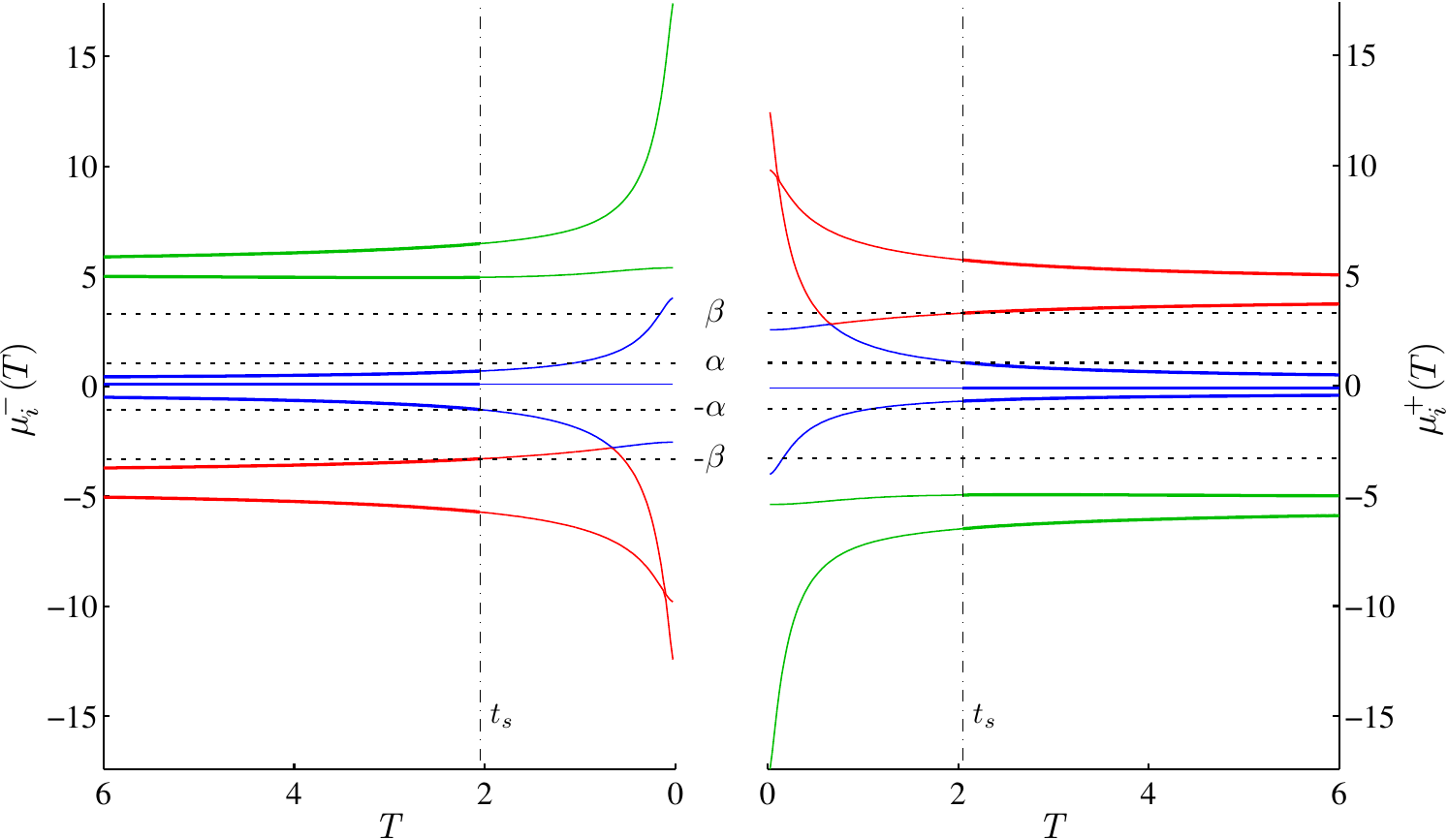}
\caption{Backward and forward FTLEs ($\mu^{\pm}_{i}$ with $i=1,...,n$). The constants $\alpha$ and
$\beta$ and the start time $t_s$ are shown.}
\label{FTLE_7DLTI}
\end{figure}

 %For this example $\hat{k}=0.7$ was chosen, leading to a start time $t_s=2.05$.

\noindent Figure \ref{FTLE_EB_7DLTI} shows the FTLEs $\mu^{\pm}_{EB}$ used to
determine the constants $\nu$, $\sigma$ as described in Section
\ref{diagnosis}. With $t_s=2.05$, for $T\in(t_s,t_c]$ we can define uniform exponential
bounds with
$\sigma=0.31$ and $\nu=3.29$. Figure \ref{FTLE_7DLTI} shows the
FTLEs $\mu^{\pm}$.

In the general case with linear-time-varying (LTV) tangent dynamics,
there is similar behavior requiring the truncation
of the time interval. The specification of the constants $t_s$ and $t_c$ can be
exclusively based on behavior of the $\mu^{\pm}$ and
$\mu_{EB}^{\pm}$; it is not necessary to determine angles within
subspaces as was done in this example to provide
insight into the root cause.

\subsection{Davis-Skodje 2D System: Attracting Slow Manifold}

Davis and Skodje (D-S) \cite{DavisSkodje} introduced a
2D nonlinear system
\begin{equation}\label{eqn:DS_ODEs}
\begin{array}{llc}
\dot{x}_1=-x_1,\\
\dot{x}_2=-\gamma x_2 + \frac{(\gamma-1)x_1+\gamma x_1^2}{(1+x_1)^2}
\end{array}
\end{equation}
defined on the state space $\{(x_1,x_2)\in \mathbb{R}^2: x_1\ge
0\hspace{2mm} {\rm and} \hspace{2mm} x_2 \ge 0\}$ with constant $\gamma>1$, which has become a benchmark for center
manifold determination. The
origin is a globally attracting equilibrium point, but more
importantly in the present context, for sufficiently large $\gamma$,
trajectories are first attracted on a faster timescale to the 1D
center manifold
\begin{equation}
\mathcal{W}^c=\{(x_1,x_2)\in \mathbb{R}^2:x_2=x_1/(1+x_1)\},
\label{slowman}
\end{equation}
 and then
follow $\mathcal{W}^c$ to the origin on a slower timescale. The two
timescales are evident in the analytic solution
\begin{equation}\label{eqn:DS_sol}
\phi(t; x_1, x_2)=\begin{bmatrix} x_1e^{-t}&\\
\\\left( x_2-\frac{x_1}{1+x_1}   \right)e^{-\gamma
t}+\frac{x_1}{1+x_1e^{-t}}e^{-t}&
\end{bmatrix}.
\end{equation}

\noindent for the flow
associated with the vector field in (\ref{eqn:DS_ODEs}). Note that if the initial state is on the center manifold, there
is no
fast timescale behavior because the coefficient of $e^{-\gamma t}$
in (\ref{eqn:DS_sol}) is zero. For this system, both the nonlinear and tangent linear dynamics have two timescales, and
it is appropriate to refer to the center manifold as the slow manifold.

The invariant center manifold $\mathcal{W}^c$ and
several other trajectories are shown in Fig.~\ref{fig:DS_traj} for
$\gamma=10$. The time interval between dots on the trajectory
is 0.1, illustrating faster motion off $\mathcal{W}^c$ than on
$\mathcal{W}^c$. From the analytical representation (\ref{slowman})
for the center manifold, we know that for any $\mathbf{x}\in
\mathcal{W}^c$,
\begin{equation} T_{\mathbf{x}}{\mathcal{W}^c}=span\{[(1+x_1)^2\quad 1]^T\}.
\end{equation}
 The linearized dynamics (\ref{lindyn}) for the D-S system have the Jacobian matrix
\begin{equation}D\mathbf{f}=\begin{bmatrix}
-1&0&\\
\frac{(\gamma-1)+(\gamma+1)x_1}{(1+x_1)^3}&-\gamma&
\end{bmatrix}.\label{eqn:ds_J}\end{equation}

\begin{figure}[ht]
\hspace{0mm}\includegraphics[scale=0.39]
{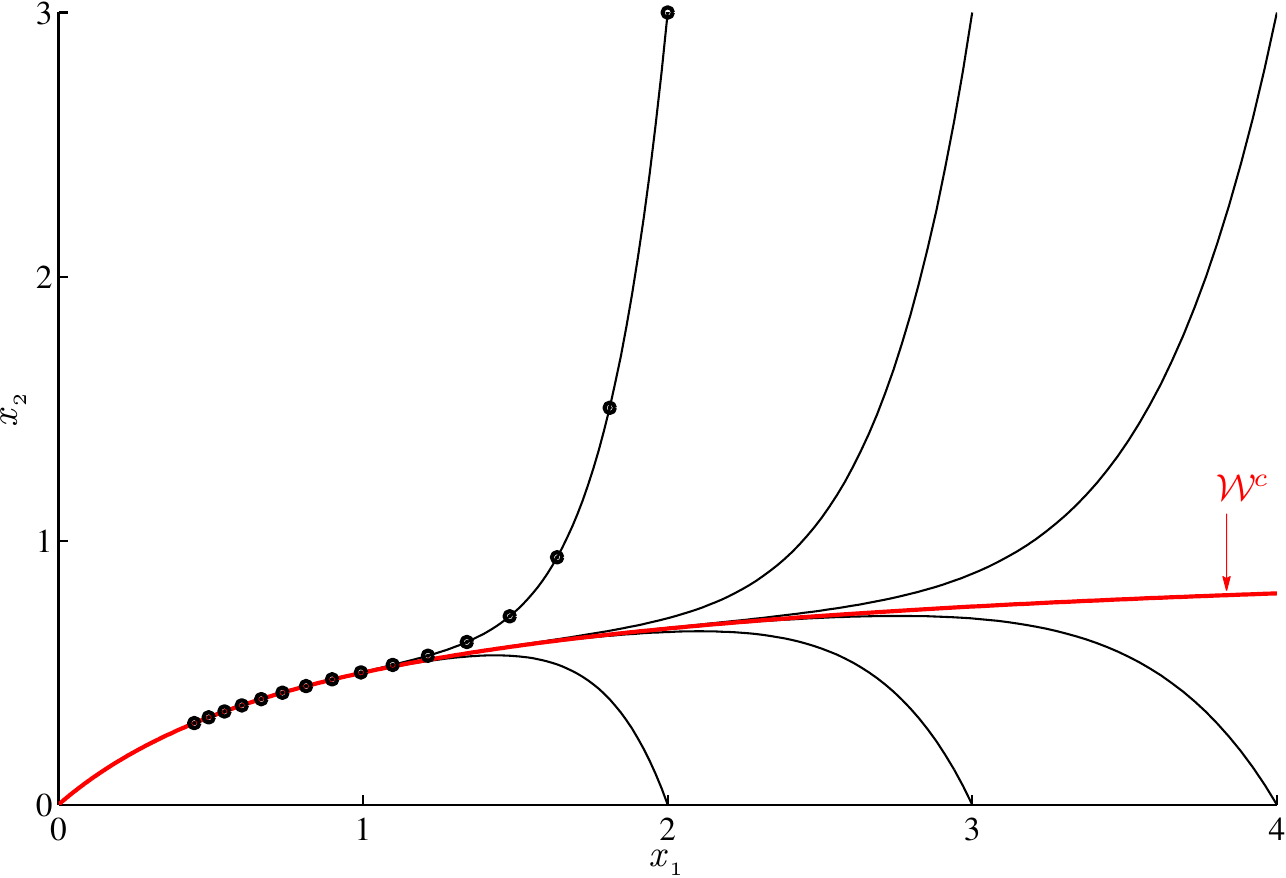}
\caption{Sample trajectories of the D-S system for $\gamma=10.0$ with the center manifold $\mathcal{W}^c$ indicated. The
dots on the trajectory departing from $\mathbf{x}=(3,2)$ are computed with $\Delta t=0.1$ and illustrate faster
motion off the center manifold than on. }
\label{fig:DS_traj}
\end{figure}

\begin{figure}[ht]%\centering
\hspace{-1.5mm}
\includegraphics[scale=.344]{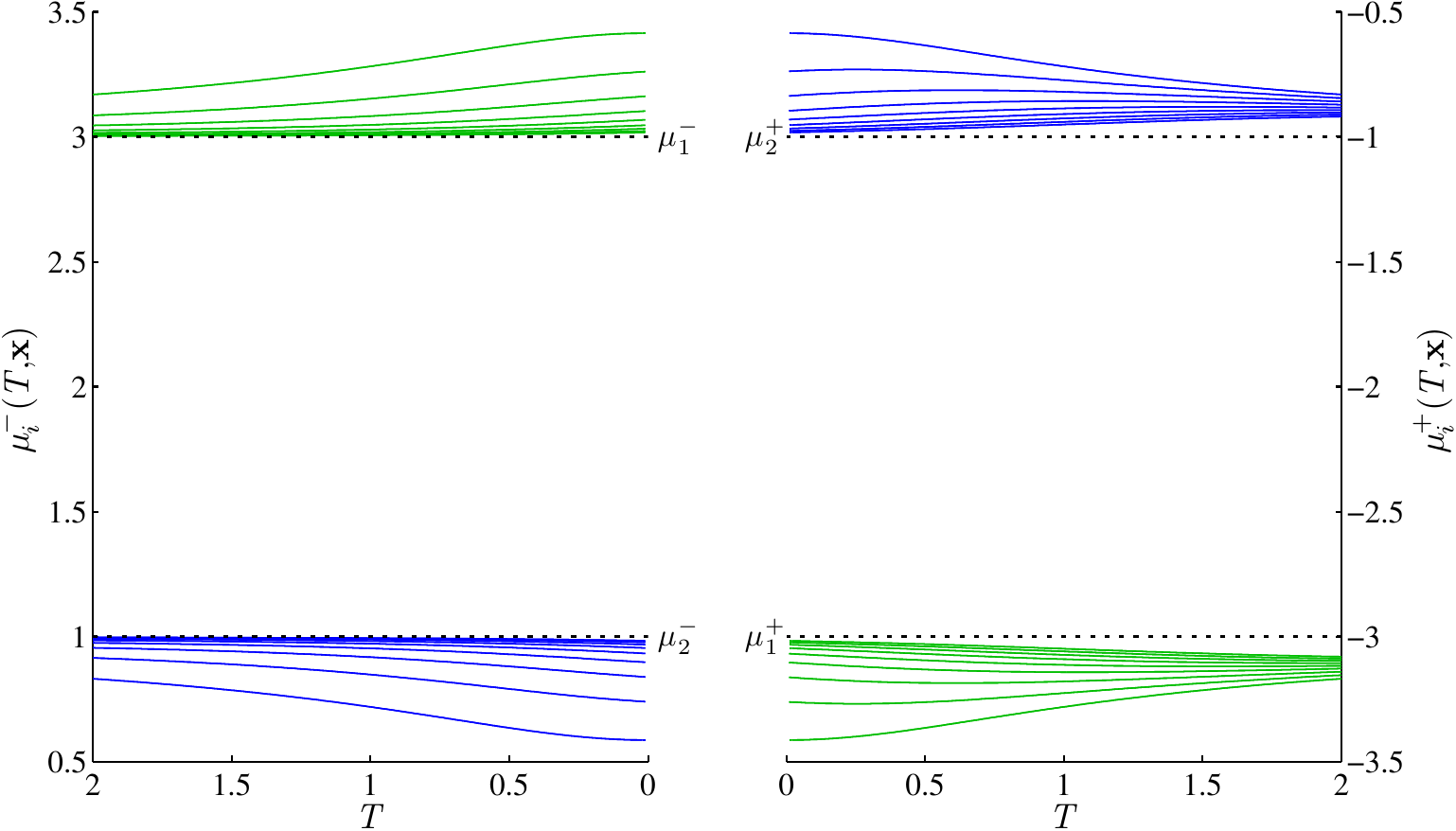}
\caption{Superposition of forward and backward FTLEs for the
Davis-Skodje system for various values of $\mathbf{x}$ illustrating uniformity. } \label{fig:DS_back
exps}
\end{figure}

Given the presence of the equilibrium point, other approaches based on eigen-analysis at the equilibrium point are
applicable: for example, integrating (\ref{eqn:DS_ODEs}) backward from an initial state perturbed slightly from the
origin in the direction of the eigenvector associated with the largest eigenvalue to compute $\mathcal{W}^c$. However
our
purpose here is to demonstrate the methodology developed in this paper, methodology that does not require the presence
of an equilibrium point.

\begin{figure}[ht]%\centering
\hspace{-4mm}
\includegraphics[scale=.36]{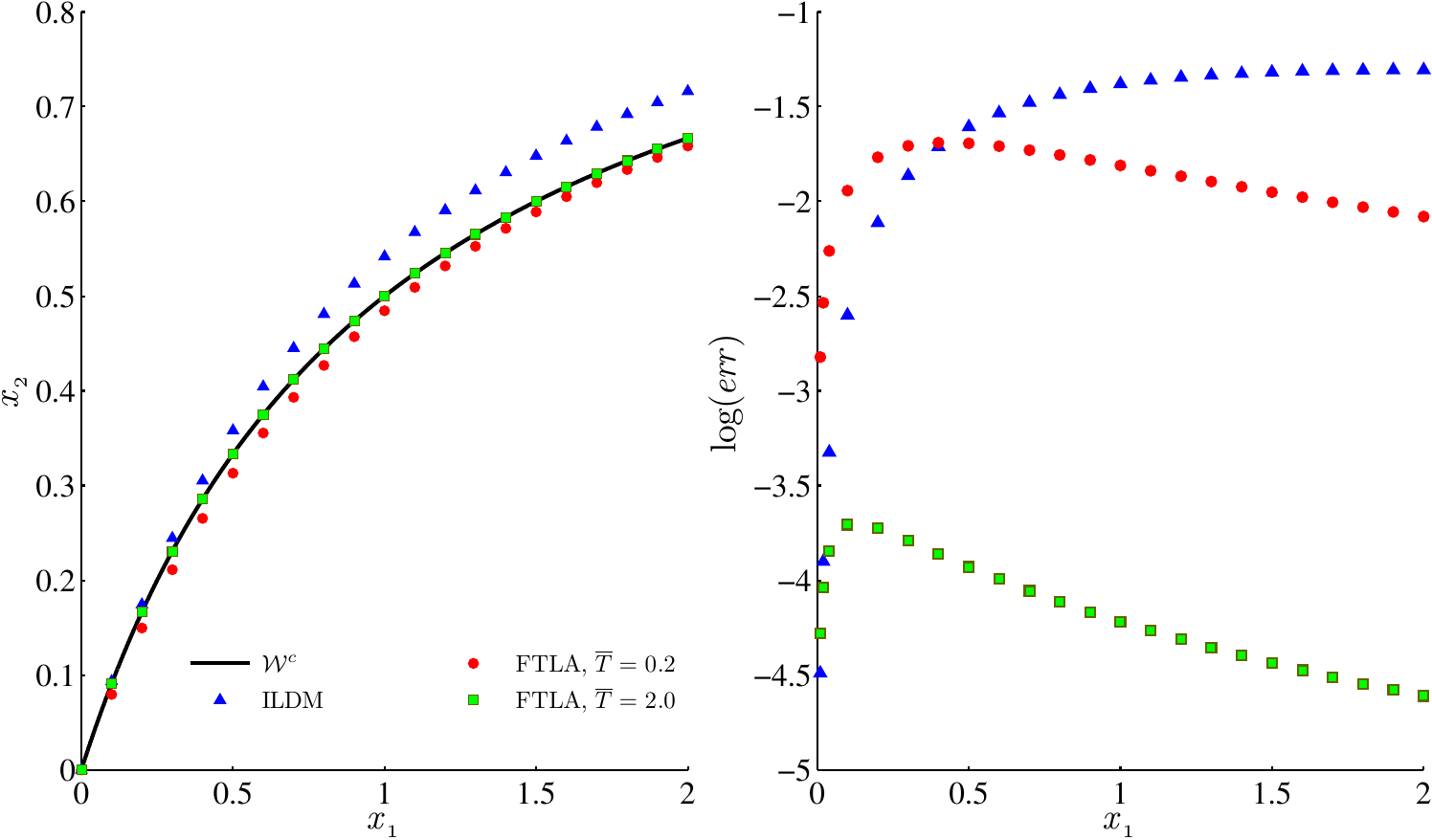}
\caption{\textit{Left plot}: Exact center invariant manifold $\mathcal{W}^c$  and its approximations calculated via FTLA
and
ILDM methods for several values of $x_1$. FTLA results are shown for two different averaging times have been used:
$\overline{T}=0.2$ and $\overline{T}=2.0$.
\textit{Right plot}: ILDM and FTLA center manifold approximation errors.}\label{SlowMan_errors_2D_DS}
\end{figure}

\subsubsection{Finite-Time Lyapunov Analysis}
 We now demonstrate the numerical application of FTLA for the case $\gamma=3$, the
case also investigated in \cite{DavisSkodje}. We consider the set
${\cal X}=\{(x_1,x_2)\in\mathbb{R}^2: 0\le x_1\le 2.0 \hspace{1.5mm} {\rm and} \hspace{1.5mm} 0\le x_2\le 1.0\}$ and
check
if the system
(\ref{eqn:DS_ODEs}), with $\gamma=3.0$, satisfies the conditions in
Def. \ref{definition1} for a finite-time uniform two-timescale
set. Figure \ref{fig:DS_back exps} shows the superposition of the
forward and backward FTLEs, as functions of $T$, for a uniform grid
of points in $\mathcal{X}$. The only possibility for two timescales is to consider $n^s=1$, $n^c=1$, $n^u=0$. Then with
$\alpha=1.0$, $\beta=3.0$, $\Delta\mu=2.0$,
$\sigma=1.0$, $\nu=3.0$,
$t_s=0$, $t_c=\overline{T}$, and $\overline{T}\ge 2.0$, the Def.
\ref{definition1} conditions are satisfied, and we conclude that
${\cal X}$ is a uniform two-timescale set resolvable over  at least $4$ convergence time constants. For the D-S system,
it
can be verified that the timescale behavior is globally uniform, so
that there is no upper limit on $\overline{T}$ unless numerical errors are an issue. The FTLVs that approximate
the fast and slow directions are $\mathbf{l}^+_1(T,\mathbf{x})$ and
$\mathbf{l}^-_2(T,\mathbf{x})$.

FTLA indicates the potential existence of a one-dimensional center manifold that can be parametrized by $x_1$. Candidate center manifold points, namely, points that satisfy
the orthogonality condition $\langle\mathbf{f}(\mathbf{x}),\mathbf{l}_1^-(\overline{T},\mathbf{x})\rangle=0$  are shown
in Fig.~\ref{SlowMan_errors_2D_DS} with $\overline{T}=0.2$ and $\overline{T}=2.0$. Over time intervals around $t_f=1$,
attraction to the center manifold occurs before the equilibrium point at the origin is reached. Because the slow and
fast timescales are not very different for $\gamma=3.0$, there is not as strong an attraction to the center manifold as
would be the case for larger values of $\gamma$, yet even for this modest level of timescale separation, the
two-timescale structure can be resolved.

\subsubsection{Asymptotic Lyapunov Analysis}
 For the D-S system, because the timescale structure is uniform on
 the entire state space,
the progress toward convergence in the first 2 units of time
continues, and it is possible to compute the asymptotic Lyapunov
exponents and vectors. The infinite-time limits of the FTLEs can be
determined analytically to be $\mu_1^+=-\gamma$, and $\mu_2^+=-1$.
The backward time limits are
$(\mu_1^-,\mu_2^-)=(\gamma,1)=(-\mu_1^+,-\mu_2^+)$.

We can analytically
compute the center FTLV $\mathbf{l}_2^-(T,\mathbf{x})$ as the eigenvector of
$\Phi(-T,\mathbf{x})^T\Phi(-T,\mathbf{x})$ corresponding to $\mu_2^-(T,\mathbf{x})$, the center
exponent in backward time. As $T$ goes to infinity,
$\mathbf{l}_2^-(T,\mathbf{x})$ can be shown to converge to

\begin{equation}\mathbf{l}_2^-(\mathbf{x})=a(x_1,x_2)
\begin{bmatrix}(1+x_1)^2\\1\end{bmatrix}
\end{equation} where $a(x_1,x_2)$ is a non-zero scalar function. For
$\mathbf{l}_2^-$ to be a unit vector, $a(x_1,x_2)$ should be
chosen appropriately. Similarly $\mathbf{l}_1^+(T, \mathbf{x})$ can be shown to converge to
\begin{equation}\mathbf{l}_1^+(\mathbf{x})=
\begin{bmatrix}0\\1\end{bmatrix}
\end{equation}
 independent of $\mathbf{x}$.

If a point $\mathbf{x}$ is on $\mathcal{W}^c$, then, using the
asymptotic Lyapunov vector $\mathbf{l}_1^-(\mathbf{x})$, the
orthogonality condition characterizing points on $\mathcal{W}^c$ is in
agreement with (\ref{slowman}). These asymptotic results lend
credence to the finite-time results, but the most important message
is that in 2 units of time, the two-timescale behavior can be
diagnosed and an accurate approximation of the center manifold can be
obtained.

\subsubsection{Invariant Center Manifold Approximation Using
Eigenvectors of $D\mathbf{f}$}

 The eigenvalues of $D\mathbf{f}$ in (\ref{eqn:ds_J}) are $-\gamma$
and $-1$; in this case they indicate the two-timescale behavior
correctly. Assuming that the span of the eigenvector, denoted
$\mathbf{e}^c$, associated with the center eigenvalue $-1$,
approximates the center subspace of the tangent plane, the ILDM method
\cite{MaasPope} estimates points on $\mathcal{W}^c$ by computing
solutions to the orthogonality condition
$\langle\mathbf{f}(\mathbf{x}),(\mathbf{e}^c)^\perp\rangle=0$. The center eigenvector $\mathbf{e}^c$
can be obtained analytically and is
\begin{equation}
\mathbf{e}^c=\begin{bmatrix} (1+x_1)^3&\\
\\1+\frac{(\gamma+1)}{(\gamma-1)}x_1&
\end{bmatrix}.
\end{equation}
The ILDM approximation to the center manifold is
\begin{equation}
x_2=\frac{x_1}{1+x_1}+\frac{2x_1^2}{\gamma^2}\left[\frac{1}{(1-\frac{1}{\gamma})(1+x_1)^3}\right].
\end{equation}

Figure \ref{SlowMan_errors_2D_DS} shows the exact manifold $\mathcal{W}^c$ along with approximations calculated with
the ILDM and FTLA methods. The ILDM
approximation is accurate around the equilibrium point (small $x_1$)
but gets worse away from the origin. The error is proportional to
$\varepsilon^2$, where $\varepsilon=1/\gamma$, consistent with the
analysis  in \cite{Kaper}. The
FTLA method provides uniformly accurate approximations when a sufficiently large averaging time is used.
$\overline{T}=2.0$ is large enough here, whereas $\overline{T}=0.2$ is not.
The center manifold approximation errors are calculated so that $err=|x_2^{\mathcal{W}^c}-\hat{x}_2|$ where
$x_2^{\mathcal{W}^c}$ is the exact $x_2$-coordinate defined in (\ref{slowman}) and $\hat{x}_2$ represents the
ILDM or FTLA $x_2$-coordinate approximation.

\subsection{3D Nonlinear System: Normally Hyperbolic Center Manifold}\label{3DnonlinEx}

Consider a nonlinear time-invariant system

\begin{equation}\label{eqn:nti_3d}
\begin{array}{ccc}
\dot{x}_1&=&ax_1,\\
\dot{x}_2&=&bx_2+\gamma(b-2a)x_1^2,\\
\dot{x}_3&=&cx_3+\gamma(c-2a)x_1^2.\\
\end{array}
\end{equation}
\noindent For the numerical results, the constants are assigned the
values $a=-0.2,\ b=-3,\ c=3$,
and $\gamma=2$.

\subsubsection{Finite-Time Lyapunov Analysis} First
the FTLEs are computed on a uniform grid on the cubic region
$\mathcal{X}=[-10, 10]^3\subset\mathbb{R}^3$. Figure \ref{fig:FTLEs} shows a superposition of all the
forward and backward FTLEs as functions of averaging time for the 36 values of $\mathbf{x}$ on the $\mathcal{X}$ grid.
The only possibility for two timescales is $n^s=n^c=n^u=1$. With $\alpha=0.8$, $\beta=3.0$, $\Delta\mu = 2.2$,
$\sigma=0.5$, $\nu=3.0$,
  $t_s=0$, $t_c=\overline{T}$ and $\overline{T}=3.0$, the Def. \ref{definition1} requirements for a uniform
two-timescale set resolvable over $6.64$ convergence time constants are satisfied.

 Having diagnosed two timescales and both fast-stable and fast-unstable
   behavior, there may be a 1D center manifold and, if so, it is normally hyperbolic.
   Because there is sufficient averaging time,
$[\mathcal{E}^c(\overline{T},\mathbf{x})]^\perp=span\{\mathbf{l}_1^-(\overline{T},\mathbf{x}),\mathbf{l}_{3}^+(\overline
{T},
\mathbf{x})\}$,
the application of the general result (\ref{prop3eq}), is a good
approximation of the orthogonal complement to the corresponding
invariant center subspace, and an accurate approximation to invariant center manifold can be obtained.

 \begin{figure}[ht]  %\centering
 \hspace{-4mm}
  \includegraphics[scale=.35] {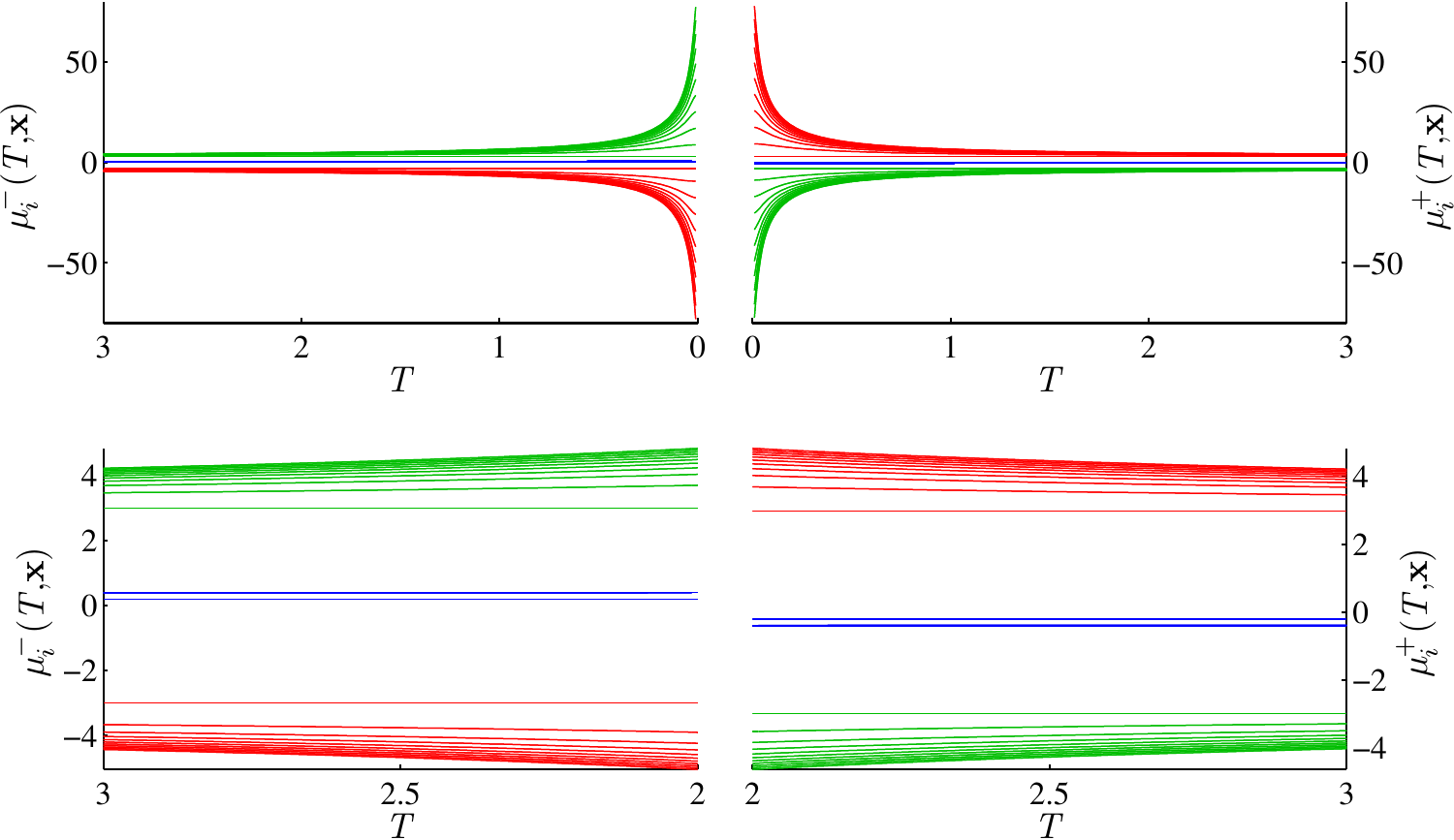}
% \mbox{\subfigure[FTLE in Forward Time]{}\quad
%        \subfigure[FTLE in Backward Time] { \includegraphics[scale=.35] {FTLE_b.pdf}}\quad}
        \caption{Forward and backward finite-time Lyapunov exponents
         for grid points on $\mathcal{X}$. For forward
        time, $\mu_1^+(T,x)$, $\mu_2^+(T,x)$, $\mu_3^+(T,x)$ are green, blue, red, resp.
For backward time, $\mu_3^-(T,x)$,
$\mu_2^-(T,x)$, $\mu_1^-(T,x)$ are red, blue, green, resp. The two lower plots zoom in on the final interval of
$T$.}\label{fig:FTLEs}
\end{figure}

After examining the FTLVs, we chose $x_1$ to parametrize $\mathcal{W}^c$, because its coordinate axis is not parallel to any of the
directions in $[\mathcal{E}^c(\overline{T},\mathbf{x})]^\perp$. For each of the values on
the grid over $x_1$, we compute the values of $x_2$ and $x_3$ that
satisfy the orthogonality conditions. {{The resulting finite-time
approximation of the postulated invariant center manifold for values of $x_1$ from -10 to
10 is plotted in Fig.~\ref{fig:Ss}.}}

\subsubsection{Invariant Center Manifold}
For this problem, there is an invariant center manifold and a means of determining it, allowing the accuracy of FTLA to
be
assessed. Over a time interval long relative to the fast timescale,
yet short relative to the slow timescale, trajectories approach the
2D manifolds $\mathcal{W}^{cu}$ and $\mathcal{W}^{cs}$, in forward and
backward time respectively, given by

\begin{equation}\label{eqn:Slow_manifolds}
\begin{array}{ccc}
\mathcal{W}^{cu}=\{(x_1,x_2,x_3)\in \mathbb{R}^3\ |\ x_2+\gamma x_1^2=0\},\\
\mathcal{W}^{cs}=\{(x_1,x_2,x_3)\in \mathbb{R}^3\ |\ x_3+\gamma x_1^2=0\}.\\
\end{array}
\end{equation}
The intersection of these sets is the invariant center manifold:
$\mathcal{W}^c$=$\mathcal{W}^{cu}\bigcap\mathcal{W}^{cs}$. These manifolds and their intersection are shown in
Fig.~\ref{fig:Ss}.

\begin{figure}[ht] % \centering
\hspace{-8mm}
\includegraphics[scale=.38] {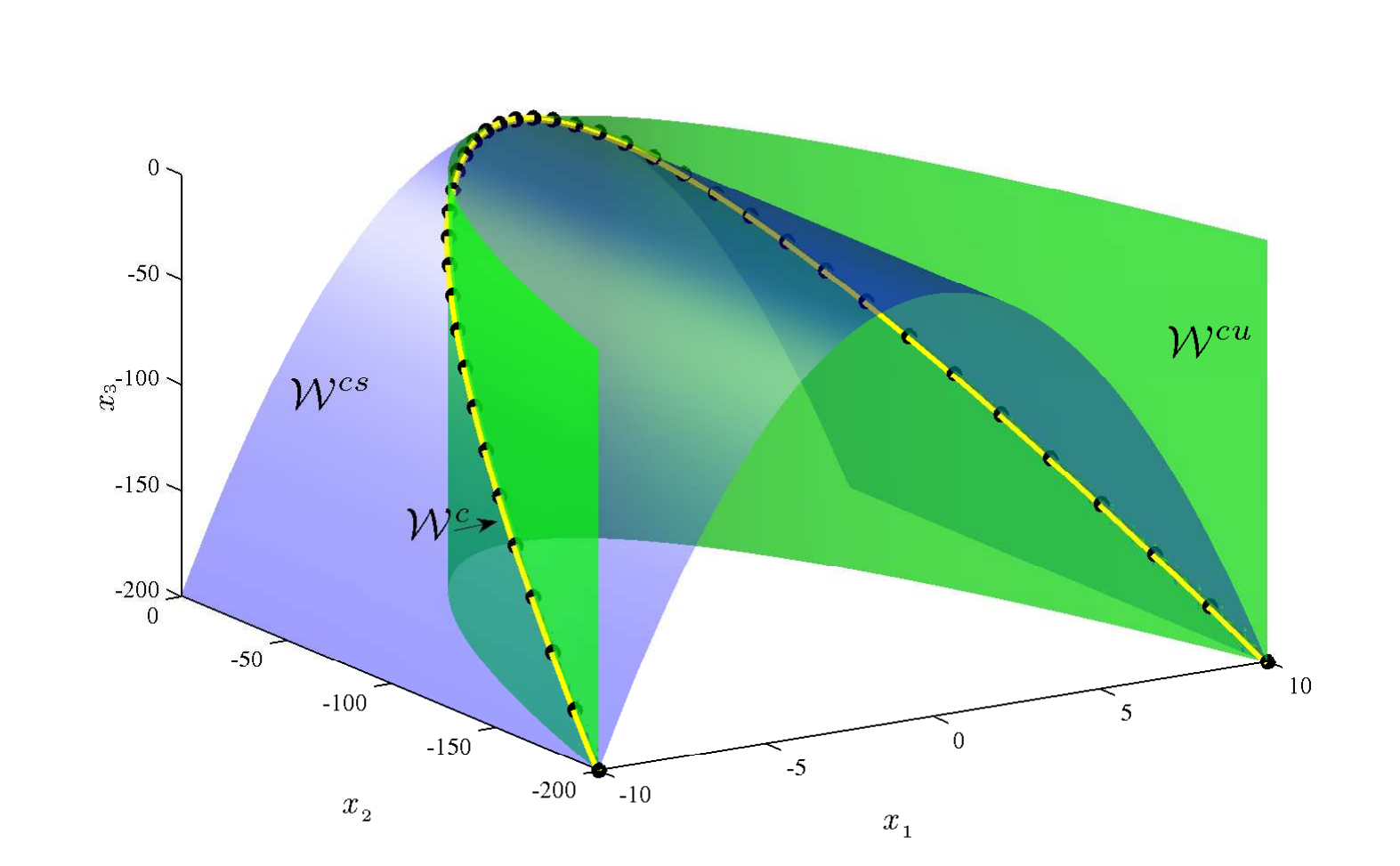}
\caption{Invariant center manifold $\mathcal{W}^c$ (yellow curve) as the intersection of the 2D manifolds
$\mathcal{W}^{cu}$ (green surface) and
$\mathcal{W}^{cs}$ (purple surface). The black rings and dots represent the FTLA approximations of points on the
invariant center manifold calculated for
$\overline{T}=3.0$.}\label{fig:Ss}
\end{figure}

At a point $\mathbf{x}\in \mathcal{W}^c$,  the vectors normal to
$\mathcal{W}^{cu}$ and $\mathcal{W}^{cs}$ are given by
$\eta_1(\mathbf{x})=[2\gamma x_1\ 1\ 0]^T$ and $\eta_2(\mathbf{x})=[2\gamma x_1\ 0\ 1]^T$
respectively. Points on $\mathcal{W}^c$, due to its invariance with
respect to the flow, satisfy the orthogonality conditions
\begin{equation}
\begin{array}{ccl}
0&=&\left\langle\eta_1(\mathbf{x}),\mathbf{f}(\mathbf{x})\right\rangle = \left\langle[2\gamma x_1\ 1\
0]^T,\mathbf{f}(\mathbf{x})\right\rangle \\
&=& 2\gamma a x_1^2+b x_2+\gamma(b-2a)x_1^2
 = b(x_2+\gamma x_1^2)\\
&&\\
0&=&\left\langle\eta_2(\mathbf{x}),\mathbf{f}(\mathbf{x})\right\rangle
 = \left\langle[2\gamma x_1\ 0\ 1]^T,\mathbf{f}(\mathbf{x})\right\rangle \\
&=& 2\gamma a x_1^2+c x_3+\gamma(c-2a)x_1^2
 = c(x_3+\gamma  x_1^2)\\
\end{array}
\end{equation}
\noindent where $\mathbf{f}(\mathbf{x})$ is the vector field given
in (\ref{eqn:nti_3d}). Figure \ref{FTLV_conv} shows $\alpha^+(T,\mathbf{x})$ and
$\alpha^-(T,\mathbf{x})$ which are respectively the angles between
$\mathbf{l}_3^+(T,\mathbf{x})$ and $\mathbf{\eta}_2(\mathbf{x})$ and between
$\mathbf{l}_1^-(T,\mathbf{x})$ and $\mathbf{\eta}_1(\mathbf{x})$. The angles are functions
of $T$ and are plotted for several values of $x_1$. As the averaging time increases, the FTLVs used to approximate
the directions of the normal vectors to the invariant center manifold align with those vectors.

For a given $x_1$, letting $(x_1,\hat{x}_2,\hat{x}_3)$ denote an
approximation of the invariant center manifold point $(x_1, -\gamma x_1^2, -\gamma x_1^2)$, we define the
approximation error to be $err=[(\hat{x}_2+\gamma x_1^2)^2+(\hat{x}_3+ \gamma x_1^2)^2]^{1/2}$.
The approximation errors for FTLA are calculated using $\overline{T}=1.0$, $\overline{T}=2.0$ and
$\overline{T}=3.0$ and plotted in Fig.~\ref{fig:errors}.

\begin{figure}
\hspace{-1.5mm}
\includegraphics[scale=.36] {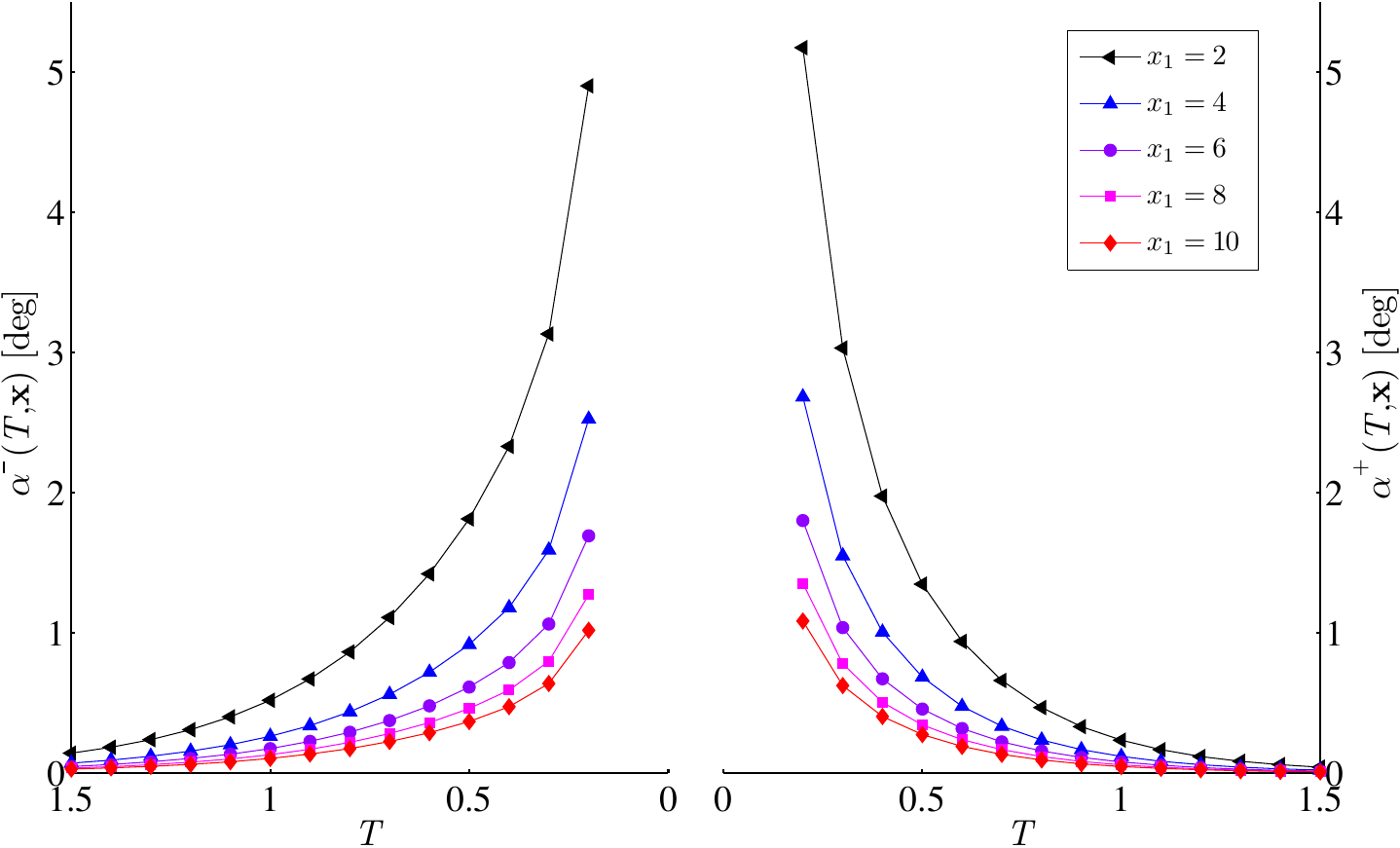}
\caption{Angles between $\mathbf{l}^-_1(T,\mathbf{x})$, $\mathbf{l}^+_3(T,\mathbf{x})$ and the directions normal to
the invariant center manifold $\mathbf{\eta}_1(\mathbf{x})$, $\mathbf{\eta}_2(\mathbf{x})$ versus the averaging time
$T$.
Points are plotted for different values of $x_1$.}\label{FTLV_conv}
\end{figure}

\begin{figure}[ht]%\centering
\hspace{-1.5mm}
\includegraphics[scale=.39] {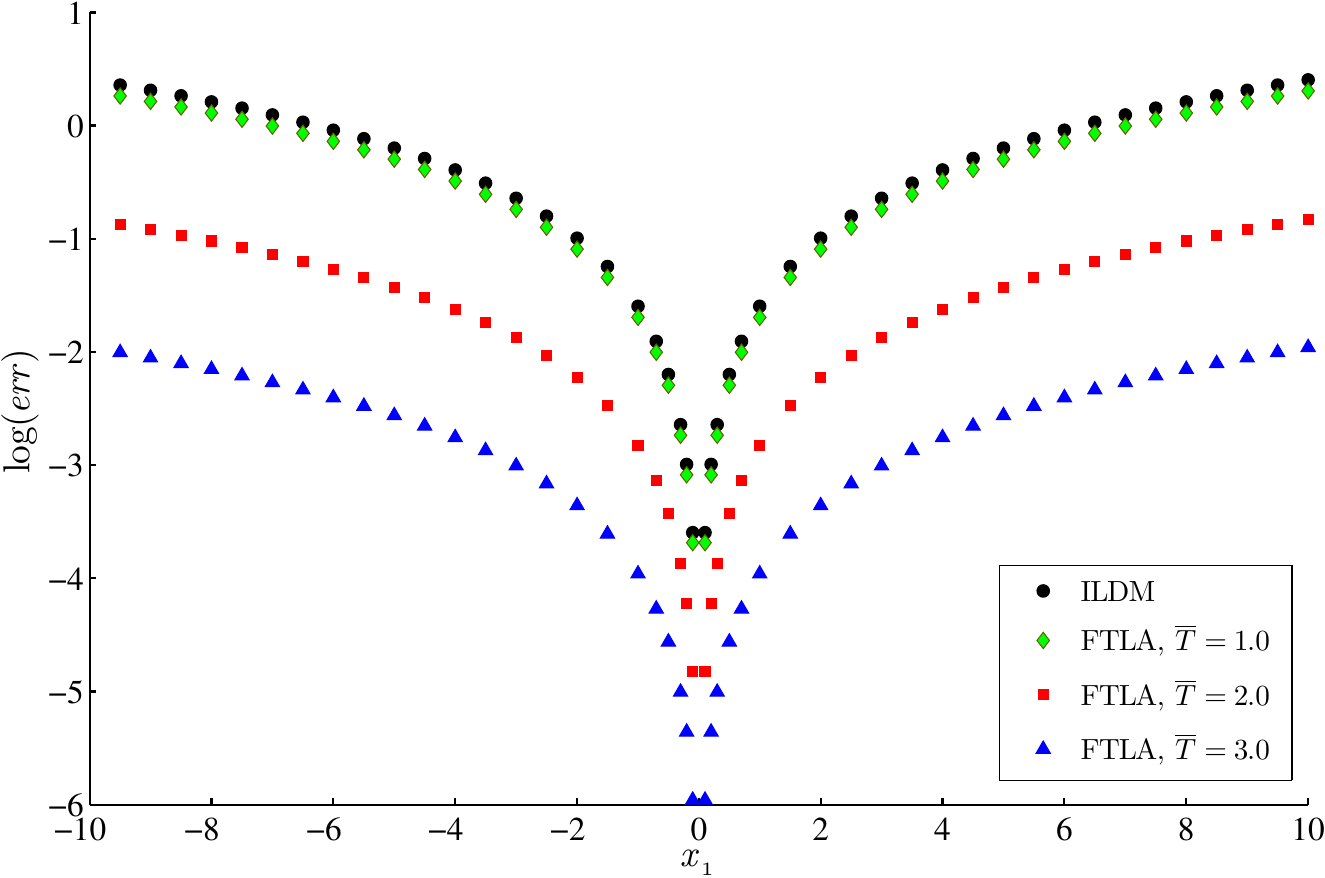}\caption{ILDM and FTLA center
manifold approximation errors calculated for various values of the independent variable $x_1$. The
FTLA approximation errors are provided for averaging times
$\overline{T}=1.0,2.0,3.0$.}\label{fig:errors}
\end{figure}

\subsubsection{Invariant Center Manifold Approximation Using
Eigenvectors of $D\mathbf{f}(\mathbf{x})$} Using the eigenvalues and eigenvectors of $D\mathbf{f}(\mathbf{x})$ (the ILDM
method \cite{MaasPope}), it
is assumed that the eigenvector corresponding to the center eigenvalue
spans the center subspace. The Jacobian matrix corresponding to the system (\ref{eqn:nti_3d}) is
\begin{equation}\label{eqn:J}D\mathbf{f}=
\begin{bmatrix} a&0&0\\2\gamma(b-2a)x_1&b&0\\
2\gamma(c-2a)x_1&0&c
\end{bmatrix}
\end{equation}
and the eigenvector corresponding to the center eigenvalue, $\lambda_c=a$ for the numerical values used,
can be written as
\begin{equation}\mathbf{v}_c=\begin{bmatrix}1,&-2\gamma x_1\left(\frac{b-2a}{b-a}\right),&
-2\gamma x_1\left(\frac{c-2a}{c-a}\right)\end{bmatrix}^T\end{equation}
Two linearly independent vectors orthogonal to $\mathbf{v}_c$ are
\begin{equation}\mathbf{w}_1=\begin{bmatrix}2\gamma x_1\left(\frac{b-2a}{b-a}\right),&1,&0\end{bmatrix}^T,
\mathbf{w}_2=\begin{bmatrix}2\gamma x_1\left(\frac{c-2a}{c-a}\right),&0,&1\end{bmatrix}^T\end{equation}
Points on the invariant center manifold are approximated using solutions to the
orthogonality conditions
\begin{equation}\label{eqn:MPA}
\langle \mathbf{w}_1,\mathbf{f}(\mathbf{x})\rangle=0,\hspace{5mm}
\langle \mathbf{w}_2,\mathbf{f}(\mathbf{x})\rangle=0
\end{equation}
 For given $x_1$, the magnitudes of the
errors in $x_2$ and $x_3$ relative to the correct values for $\mathcal{W}^c$ are
$2\gamma x_1^2 \frac{a^2}{(a-b)b}$ and
$2\gamma x_1^2 \frac{a^2}{(a-c)c}$ respectively. Taking the norm of these errors, the center
manifold approximation error for the
ILDM method is plotted in Fig.~\ref{fig:errors}. The ILDM error is
similar to that for FTLA when the averaging time is
$\overline{T}=1.0$, but FTLA gives greater accuracy for the longer averaging
times  $\overline{T}=2.0$ and $\overline{T}=3.0$.

%%%%%%%%%%%%%%%    MCK   %%%%%%%%%%%%%%%%%%%%%
\subsection{4D Hamiltonian System: Mass-Spring-Damper System}\label{MCKexample}

To demonstrate the use of FTLA to locate points on a two-dimensional normally hyperbolic center manifold,
we consider the optimal control of a mass-(nonlinear) spring-damper system modeled as
\begin{equation}\label{eqn:MCK5_model}\begin{array}{ccl}
\dot{x}_1&=&x_2,\\
\dot{x}_2&=&-\frac{1}{m}\left(c x_2+k_1x_1+k_2x_1^3\right)+\frac{u}{m} ,\\
\end{array}\end{equation}
where $x_1$ is the displacement of the mass $m$ measured from the rest position of the spring, $u$ is the applied scalar
control,
$k_1$ and $k_2$ are the coefficients of the linear and cubic contributions to the spring force, and $c$ is the
damping
coefficient. For the problem of minimizing the function
\begin{equation}\label{eqn:MCK6_model}\begin{array}{ccl}
min\:\:\:J=\int_0^{t_f}\frac{1}{2}u^2\, \mathrm{d}t\,,\\
\end{array}\end{equation}
subject to the dynamic constraint (\ref{eqn:MCK5_model}) and specified initial and final conditions on $x_1$
at a specified final time $t_f$, Pontryagin's minimum principle leads to first-order necessary conditions in the form of
a
boundary value problem for the Hamiltonian system
\begin{equation}\label{eqn:MCK9_model}\begin{array}{ccl}
\dot{x}_1&=&x_2 ,\\
\dot{x}_2&=&-\frac{1}{m}\left(cx_2+k_1x_1+k_2x_1^3+\frac{\lambda_2}{m}\right) ,\\
\dot{\lambda}_1&=&\frac{\lambda_2}{m}\left(k_1+3k_2x_1^2\right) ,\\
\dot{\lambda}_2&=&-\lambda_1+c\frac{\lambda_2}{m} ,\\
\end{array}\end{equation}
where ${\lambda}_1$ and ${\lambda}_2$ are adjoint variables and the minimizing control is $u^*=-\lambda_2/m$.
For consistency with the rest of the paper, we consider (\ref{eqn:MCK9_model}) in the form $\dot{\mathbf{x}} =
\mathbf{f}(\mathbf{x})$ with $\mathbf{x}=[x_1\,,\,x_2\,,\,\lambda_1\,,\,\lambda_2]^T\in\mathbb{R}^4$ and $\mathbf{f}$
defined appropriately.

For small values of $m$, the Hamiltonian system is in singularly perturbed standard form \cite{petar}, and the system
can be expected to evolve on disparate timescales. Using the two-timescale geometry to solve the boundary-value problem
has been addressed in
\cite{ACC09,Guck09,Rao99,Kopell}. Here we focus on applying FTLA to the Hamiltonian system
(\ref{eqn:MCK9_model}) to diagnose two-timescale behavior and locate points on the center manifold, which is in this
case a slow manifold. The linearized dynamics (\ref{lindyn}) have the Jacobian matrix
\begin{equation}\label{eqn:MCK20_model}\begin{array}{ccl}
D\mathbf{f}&=&\begin{bmatrix}
               0&1&0&0& \\
\frac{1}{m}\left(-k_1-3k_2x_1^2\right)&-\frac{c}{m}&0&-\left(\frac{1}{m}\right)^2 \\
\frac{\lambda_2}{m}\left(6k_2x_1\right)&0&0&\frac{1}{m}\left(k_1+3k_2x_1^2\right) \\
               0&0&-1&\frac{c}{m}
              \end{bmatrix}\ .\\
\end{array}\end{equation}
For the numerical results we use $m=0.5\:,\:k_1=1\:,\:k_2=0.01\:$, and $\:c=4\sqrt{k_1m}$.

\subsubsection{Finite-Time Lyapunov Analysis}
FTLA is applied in a region $\mathcal{X}=(-1.0,6.0)\times(-5.0,-1.9)\times(7.0,15.0)\times(0.8,5.0)$, chosen such that
the ILDM method is applicable
(i.e., the eigenvalues of $D\mathbf{f}$ are real),
yet the center manifold curvature is large enough that the ILDM method produces noticeable errors.
We present results for the five points: $\mathbf{x}_{1}=[3.00\:,\:-2.0\:,\:7.5\:,\:2.0]^T$,
$\mathbf{x}_{2}=[2.85\:,\:-2.0\:,\:9.3\:,\:2.0]^T$,
$\mathbf{x}_{3}=[2.70\:,\:-2.0\:,\:11.0\:,\:2.0]^T$, $\mathbf{x}_{4}=[2.55\:,\:-2.0\:,\:12.8\:,\:2.0]^T$,
and $\mathbf{x}_{5}=[2.40\:,\:-2.0\:,\:14.5\:,\:2.0]^T,$
that are representative of all the points in $\mathcal{X}$.
Figure \ref{fig:MCK_FTLE} shows the forward and backward Lyapunov exponents for the five points
as functions of the averaging time $T$. Because the system is Hamiltonian, the FTLEs should be symmetric about the
origin. With $n^s=n^u=1$, $n^c=2$, $\alpha=0.52$, $\beta=5.64$, $\Delta\mu = 5.12$, $\sigma=0.66$, $\nu=5.19$, $t_s=0$
and
$t_c=\overline{T}=0.50$,
the conditions given in Def. \ref{definition1} for a uniform two-timescale set resolvable over $2.6$ convergence
time constants are satisfied. Figure~\ref{fig:MCK_EXP_BNDS} shows the FTLEs and exponential bounds that
were computed as described in Section \ref{diagnosis}.

\begin{figure}[tbh]%\centering
\hspace{-1.5mm}
\includegraphics[scale=.35]{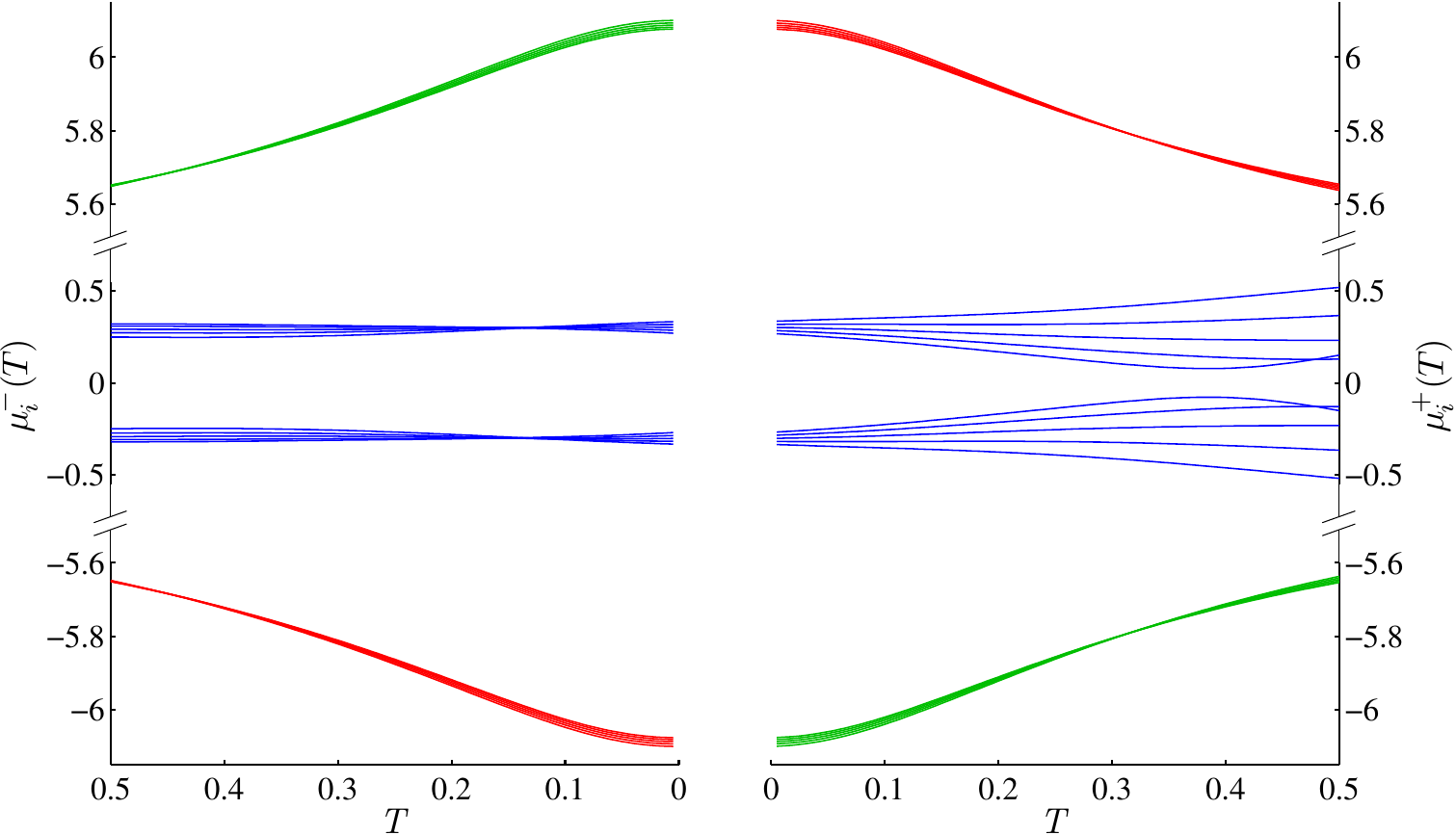}
\caption{Superposition of backward and forward FTLEs
for points $\mathbf{x}_1,\mathbf{x}_2,\mathbf{x}_3,\mathbf{x}_4,$ and $\mathbf{x}_5$. Note that only segments of the
y-axis are shown to highlight the center FTLEs.}\label{fig:MCK_FTLE}
\end{figure}

\begin{figure}[tbh]%\centering
\hspace{-1.7mm}
\includegraphics[scale=.348]{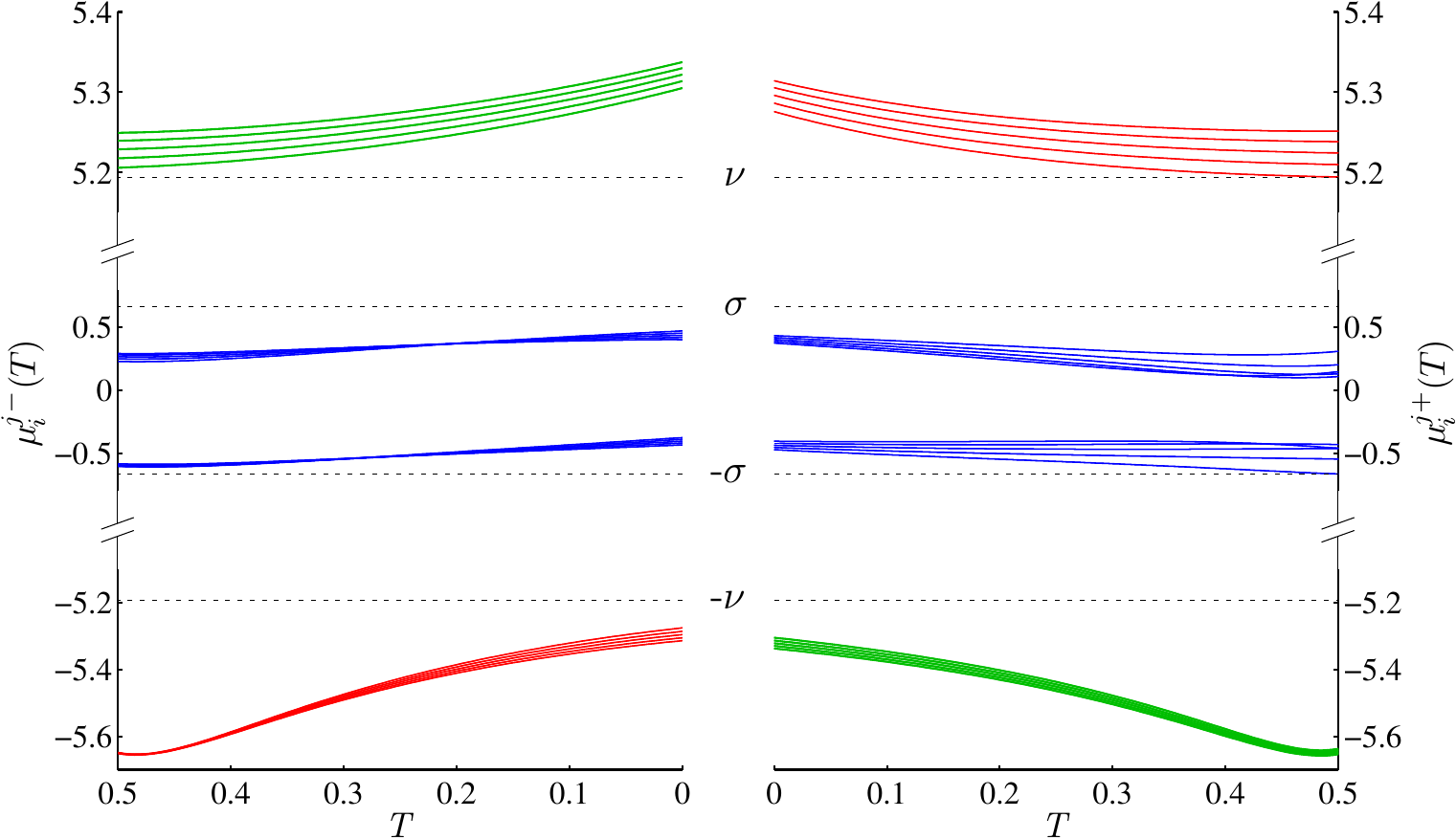}
\caption{FTLEs $\mu^{j \pm}_i$ with $i=1,...,n^j$ and $j=s,c,u$ for the subspaces $\mathcal{E}^s(0.5,\mathbf{x})$,
$\mathcal{E}^c(0.5,\mathbf{x})$, $\mathcal{E}^u(0.5,\mathbf{x})$ and determination of the constants $\nu$ and $\sigma$
for
$\mathbf{x}_1$, $\mathbf{x}_2$, $\mathbf{x}_3$, $\mathbf{x}_4$, and $\mathbf{x}_5$ as functions of
time. The distance between $\nu$ and $\sigma$ is actually larger than it appears since only segments of the vertical
axis are
shown.}\label{fig:MCK_EXP_BNDS}
\end{figure}

\subsubsection{Computing Center Manifold Points Using FTLA}\label{mckManifold}
The center subspace $\mathcal{E}^c(\overline{T},\mathbf{x})$ has dimension $n^c=2$ and can be written as (see
(\ref{EdefsT}))
\begin{equation}\label{eqn:MCK10_model}
\begin{array}{lcl}
\mathcal{E}^c(\overline{T},\mathbf{x})&=&{\cal L}^+_3(\overline{T},\mathbf{x})\cap {\cal
L}^-_2(\overline{T},\mathbf{x})\\
\end{array} \end{equation} with its orthogonal complement (\ref{prop3eq}) given by
\begin{equation}\label{eqn:MCK11_model}
[\mathcal{E}^c(\overline{T},\mathbf{x})]^\bot
=span\{\mathbf{l}^-_1(\overline{T},\mathbf{x}),\mathbf{l}^+_4(\overline{T},\mathbf{x})\}.
\end{equation}
The existence of a 2D center manifold is postulated. As described in Section \ref{theory}, $n^c$ coordinates are chosen to parametrize $\mathcal{W}^c$ such that their coordinates
axes are not parallel to any of the directions in $[\mathcal{E}^c(\overline{T},\mathbf{x})]^\perp$, namely in
$\mathbf{l}^-_1(\overline{T},\mathbf{x})$ and $\mathbf{l}^+_4(\overline{T},\mathbf{x})$ directions. For example
\begin{equation}
\begin{array}{clcl}
\mathbf{l}_1^-(0.5,\mathbf{x}_1)=&[0.33\;,\;0.89\;,\;0.05\;,\;0.31]^T,\\
\mathbf{l}_4^+(0.5,\mathbf{x}_1)=&[-0.01\;,\;0.00\;,\;-0.16\;,\;0.99]^T. \\
\end{array}
\end{equation}
The directions of $x_2$ and $\lambda_2$ are almost parallel respectively to $\mathbf{l}_1^-$ and $\mathbf{l}_4^+$,
so we choose the independent variables to be $x_1$ and $\lambda_1$.
We use the $(x_1,\lambda_1)$ coordinates of the five points $\mathbf{x}_j, j=1,\dots,5$ as the grid in the independent
coordinate plane and compute the  $(x_2, \lambda_2)$ coordinates for the graph of $\mathcal{W}^c(T)$ by solving the
orthogonality conditions.

For Def.~\ref{definition1}, the value of $\overline{T}$ must apply at each point in $\mathcal{X}$; to do so, it must be
the minimum over all the maximum forward and backward averaging times on $\mathcal{X}$. It can be
beneficial in computing center manifold points to use averaging times greater than  $\overline{T}$ when possible. An iterative procedure for determining the averaging time during convergence toward the center manifold is described in \ref{AppendixB}. For the converged points, the forward and backward averaging times were increased to 5.0 and 2.0 respectively.

Because the exact location of the invariant center manifold is not known, we use the following means to assess accuracy.
The estimated invariant center manifold points $\hat{\mathbf{x}}_j, j=1,\dots,5$ are propagated backward and
forward in time to $\phi(t^{\pm},\hat{\mathbf{x}}_j)$. Then for each of the end points, we fix the independent
variables, $x_1$ and $\lambda_1$, and use FTLA to recompute the dependent variables, $x_2$ and $\lambda_2$
for the center manifold point estimate. If the FTLA method computed points on the invariant center manifold without
error, then the propagated estimates and re-estimated points would be the same; the degree of inconsistency is thus an
indication of accuracy and invariance. The same procedure is performed for the ILDM estimates.

Figure \ref{traj_l1_x2}, showing points and trajectories projected onto the $\lambda_1$-$x_2$ plane, indicates that FTLA
is much more
consistent than the ILDM method. The trajectories departing from initial points calculated with FTLA (black circles)
propagate to points (black squares forward and black diamonds backward) close to those re-estimated, the interpretation
being that by starting closer to the invariant center manifold the trajectories follow the center manifold for a longer
time.
Although the initial ILDM points (red circles) appear close to the initial FTLA points, the high degree of inconsistency
at
the end points indicates greater inaccuracy.

\begin{figure}[tbh]%\centering
\hspace{-1mm}
\includegraphics[scale=.40]{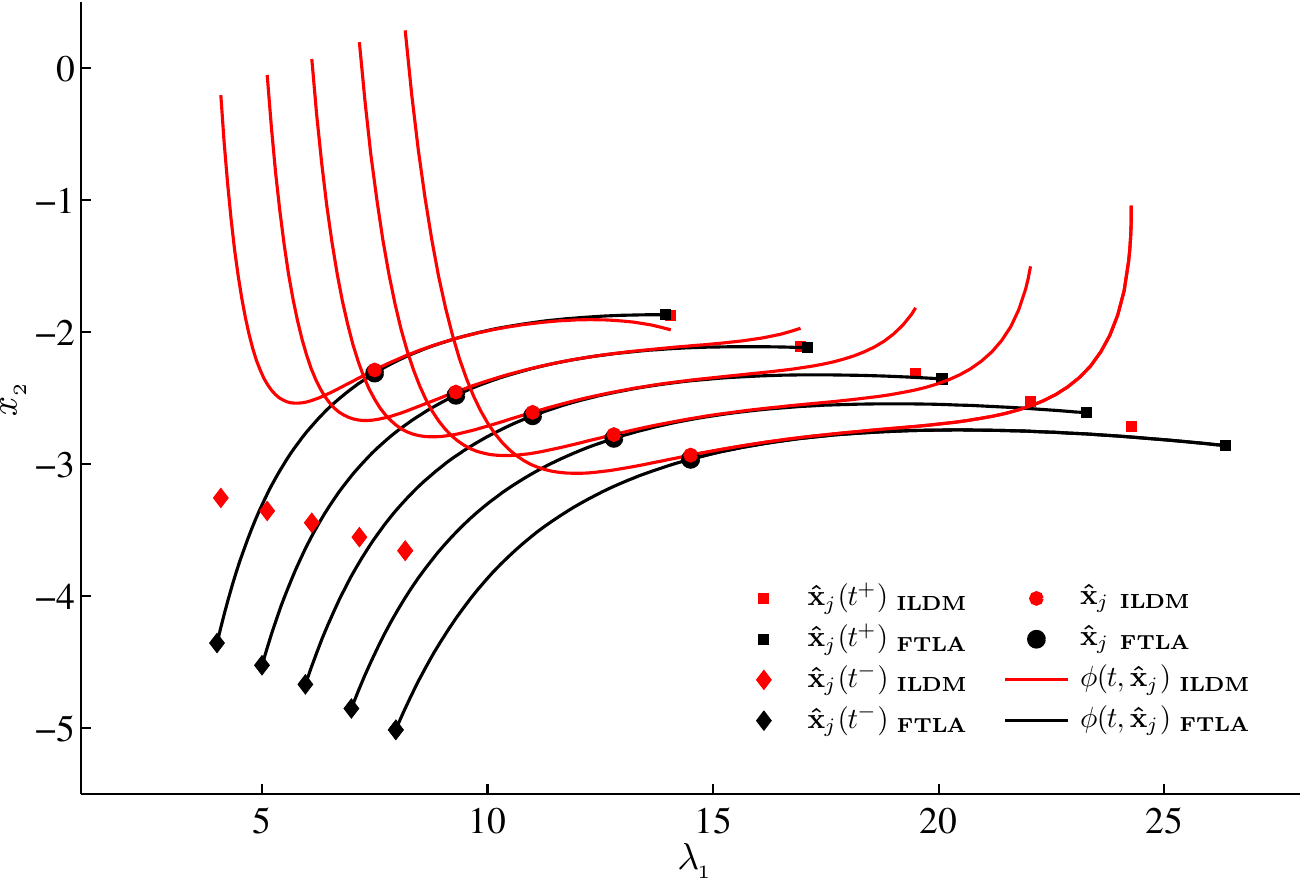}
\caption{Projection onto the $\lambda_1$-$x_2$ plane of the forward and backward propagations from initial points on
the center manifold (circles). The independent coordinates of the points at the end of the trajectories are used to
compute new estimates on the center manifold (diamonds-backward, squares-forward). Points in black
refer to estimates calculated via FTLA while the lighter ones are computed with ILDM.}\label{traj_l1_x2}
\end{figure}

Table~\ref{IP} shows quantitatively the center manifold estimation error for the FTLA and ILDM methods.
An invariance error percent (IP) is defined by
\begin{equation}
\label{eqn:MCK15_model}
IP_{x_2}^{\pm} = \frac{\|x_2(\hat{\mathbf{x}}_j(t^\pm))-x_2(\phi(t^{\pm},\hat{\mathbf{x}}_j))\|}{
\|x_2(\hat{\mathbf{x}}_j(t^\pm))\|}*100 \\
\end{equation}

\noindent where $\hat{\mathbf{x}}_j(t^+)$ (squares) and $\hat{\mathbf{x}}_j(t^-)$ (diamonds)
are estimates of points on the center manifold calculated respectively from $\phi(t^{+},\hat{\mathbf{x}}_j)$ and
$\phi(t^{-},\hat{\mathbf{x}}_j)$ via FTLA or ILDM. The trajectory end points
$\phi(t^{\pm},\hat{\mathbf{x}}_j)$  in Fig.~\ref{traj_l1_x2}) are for $t^+=1.5$ and $t^-=-1.0$). Finally $x_2(\cdot)$
denotes the $x_2$ coordinate of argument.
The explanation for $IP_{\lambda_2}^{\pm}$ is analogous. The $IP$ values indicate that FTLA produces accurate
approximations to points on the invariant center manifold
and is significantly more accurate than the ILDM method.

\begin{table} \caption{Invariance error percent for $x_2$ and $\lambda_2$ for FTLA and ILDM methods.}
\label{IP}
\begin{adjustwidth}{-0.21cm}{}
  \begin{tabular}{c c c c c c c c c}
    & \multicolumn{2}{c}{\scriptsize{$IP_{x_2}^+$}} & \multicolumn{2}{c}{\scriptsize{$IP_{x_2}^-$}} &
\multicolumn{2}{c}{\scriptsize{$IP_{\lambda_2}^+$}} & \multicolumn{2}{c}{\scriptsize{$IP_{\lambda_2}^-$}}
\\  \cmidrule(lr){2-9} %\cline{2-9}
    & \tiny{$FTLA$} & \tiny{$ILDM$} & \tiny{$FTLA$} & \tiny{$ILDM$} &
\tiny{$FTLA$} & \tiny{$ILDM$} & \tiny{$FTLA$} & \tiny{$ILDM$}  \\
    \cmidrule(lr){1-9}
    \multicolumn{1}{c}{\scriptsize {$\hat{\mathbf{x}}_1$}}& \scriptsize {1E-4} & \scriptsize {5.7E0} & \scriptsize
{5.1E-2} &
\scriptsize {9.4E1} & \scriptsize {3.0E-4} &  \scriptsize {1.1E1} & \scriptsize {6.0E-4} & \scriptsize {5.0E-1}  \\
    \multicolumn{1}{c}{\scriptsize {$\hat{\mathbf{x}}_2$}}& \scriptsize {3.0E-4} & \scriptsize {6.5E0} & \scriptsize
{2.5E-2} &
\scriptsize {9.8E1} & \scriptsize {5.0E-4} & \scriptsize {1.2E1} & \scriptsize {3.0E-4} & \scriptsize {5.6E-1}  \\
    \multicolumn{1}{c}{\scriptsize {$\hat{\mathbf{x}}_3$}}& \scriptsize {1.0E-4} & \scriptsize {2.2E1} & \scriptsize
{2.0E-3} &
\scriptsize {1.0E2} & \scriptsize {2E-4} & \scriptsize {3.8E1} & \scriptsize {$<$1E-5} & \scriptsize {6.0E-1}  \\
    \multicolumn{1}{c}{\scriptsize {$\hat{\mathbf{x}}_4$}}& \scriptsize {2.8E-3} & \scriptsize {4.1E1} & \scriptsize
{3.5E-2} &
\scriptsize {1.1E2} & \scriptsize {4.6E-3} & \scriptsize {6.9E1} & \scriptsize {4.0E-4} & \scriptsize {6.5E-1}  \\
    \multicolumn{1}{c}{\scriptsize {$\hat{\mathbf{x}}_5$}}& \scriptsize {1.4E-2} & \scriptsize {6.2E1} & \scriptsize
{8.0E-2} &
\scriptsize {1.1E2} & \scriptsize {2.2E-2} & \scriptsize {1.0E2} & \scriptsize {1.0E-4} & \scriptsize {6.9E-1}  \\
\cmidrule(lr){1-9}
\end{tabular}
\end{adjustwidth}
\end{table}

\section{Conclusions} \label{conclusions}
The practical goal of this work was to use finite-time Lyapunov analysis to improve accuracy, and extend applicability, relative to the intrinsic low-dimensional manifold method,
in estimating points on slow manifolds and more generally normally hyperbolic center manifolds. This has been accomplished as demonstrated in several examples of increasing dimension and complexity. In
addition, a definition of a uniform finite-time two-timescale set has been
proposed with requirements on the finite-time Lyapunov spectrum and the subspaces constructed from the finite-time
Lyapunov vectors, accounting for finite-time features -
non-modal growth and rate of subspace convergence to the desired invariant subspaces. Although the examples show there exist systems for which the finite-time Lyapunov analysis method is viable, further experience is needed to clarify how broadly applicable it is.

 Finite-time Lyapunov analysis of
the tangent linear dynamics provides an alternative diagnostic
approach to eigen-analysis of the associated system matrix (the
Jacobian matrix associated with the vector field). Though we have
used this finite-time information for approximating
points on invariant manifolds, the finite-time information could
potentially be used (a) to suggest a transformation of coordinates
leading to the standard form required for the analytical singular
perturbation approach, or more generally to coordinates adapted to the manifold structure, (b) to guide the
selection of independent and dependent variables in the application
of the quasi-steady-state approximation, the zero-derivative approach, and the Roussel-Fraser partial differential equation approach, and (c) to obtain an
invariant manifold approximation that could subsequently be refined by another method. Also in the solution of
boundary-value problems for two-timescale systems, determining
points on manifolds to approximate certain missing boundary conditions at each end is exactly what is needed.

\section*{Acknowledgments} Stimulating discussions
with S.-H. Lam started the first author on this research. Discussions with L.-S. Young and A. Gorodetski were
instrumental for understanding the relevant dynamical systems theory. Helpful
discussions with Y. B. Pesin and B. Villac are also
acknowledged.


\appendix
\section{Subspace Convergence}
\label{appendixA}
Proposition \ref{theorem1} below gives the
exponential rate at which the finite-time Lyapunov subspaces,
introduced in Section \ref{FTLE} and expressed in terms of the FTLVs, evolve
with increasing $T$ toward their asymptotic limits, under hypotheses in which these  limits exist. Most of the
ideas in Proposition \ref{theorem1} and its proof can be found in
\cite{Ershov,Orszag}. The new element here is that convergence of a
particular Lyapunov subspace is addressed explicitly, rather than
the convergence of individual Lyapunov vectors (see \cite{Haller11} for an alternative approach for a special case of a
co-dimension one subspace).

\begin{definition}\label{strongly non-degen} \cite{Orszag} The Lyapunov spectrum is
 strongly non-degenerate at a point $\mathbf{x}$, if there exists positive constants
$t_s$ and $\delta$ such that the spectral gap between each
neighboring pair of forward FTLEs,
$\mu^+_{i+1}(T,\mathbf{x})-\mu^+_{i}(T,\mathbf{x})$,
$i=1,\dots,n-1$, is greater than $\delta$ for all $T
> t_s$ and likewise for the backward exponents.\end{definition}

To consider the convergence of a Lyapunov subspace
$\mathcal{L}^+_j(T,\mathbf{x})$ with $T$, we focus on a particular
spectral gap and bound it for use in the proposition
that follows.
\begin{definition}\label{rel gap}[Spectral Gap Lower Bound] For a specified $t_s > 0$, the lower bound on the
spectral gap ${\Delta\mu}^+_j(\mathbf{x})$ between
neighboring forward FTLEs $\mu^+_{j}{(T,\mathbf{x})}$ and
$\mu^+_{j+1}{(T,\mathbf{x})}$, for a particular
$j\in\{1,2,\dots,n-1\}$, is
\begin{equation} {\Delta\mu}^+_j(\mathbf{x}) :=
\inf_{T\ge t_s}(\mu^+_{j+1}(T,\mathbf{x})-\mu^+_{j}(T,\mathbf{x})).
\end{equation}
Similarly the spectral gap bound
${\Delta\mu}^-_k(\mathbf{x})$ between neighboring backward
FTLEs $\mu^-_{k-1}{(T,\mathbf{x})}$ and $\mu^-_{k}{(T,\mathbf{x})}$
is defined as
\begin{equation} {\Delta\mu}^-_k(\mathbf{x}):=\inf_{T>t_s}(\mu^-_{k-1}(T,\mathbf{x})-\mu^-_{k}(T,\mathbf{x})).
\end{equation}\end{definition} \noindent

\begin{proposition}\label{theorem1} Consider the dynamical system (\ref{nldyn}) on a compact invariant
subset $\mathcal{Y}$ of the state space $\mathbb{R}^n$. At a
Lyapunov regular point $\mathbf{x}\in \mathcal{Y}$ for which there
exists $t_s>0$ and $\delta>0$ such that the Lyapunov spectrum is
strongly non-degenerate for $T>t_s$ and for which there is a
nonzero lower bound ${\Delta\mu}^+_j(\mathbf{x})$ on the spectral gap for a
specific value of $j$, the subspace $ {\cal L}^+_j(T,\mathbf{x})$
approaches the fixed subspace ${\cal L}^+_j(\mathbf{x})$, defined in
Section \ref{asympt} in terms of the asymptotic Lyapunov exponent
$\mu^+_j(\mathbf{x})$. It approaches at an exponential rate characterized, for
every sufficiently small $\Delta T
>0$, by
\begin{equation}\label{expbdd_A} \textrm{dist}({\cal L}_j^+(T,\mathbf{x}),
{\cal L}_j^+(T+\Delta T,\mathbf{x}))\le K
e^{-{\Delta\mu}^+_j(\mathbf{x})\cdot T},\end{equation} for
all $T>t_s$, where $K>0$ is $\Delta T$ dependent but $T$
independent. Similarly, as $T$ increases, the subspace $ {\cal
L}^-_k(T,\mathbf{x})$ approaches the fixed subspace ${\cal
L}^-_k(\mathbf{x})$ at a rate proportional to
$exp(-\Delta\mu^-_k(\mathbf{x}) \cdot T)$. \end{proposition}

{\it Proof of Proposition \ref{theorem1}:} Using (\ref{distance}) we have

\begin{equation}\label{inners}
\begin{aligned}
&\textrm{dist}({\cal L}_j^+(T,\mathbf{x}),{\cal L}_j^+(T+\Delta
 T,\mathbf{x}))=\|L_j^+(T,\mathbf{x})^T L_{j'}^+(T+\Delta T,\mathbf{x})\|_2\\ &=\left\|\left[%
	\begin{array}{c}
	    \mathbf{l}_1^+(T,\mathbf{x})^T \\
	    \mathbf{l}_2^+(T,\mathbf{x})^T \\
	    \vdots \\
	    \mathbf{l}_j^+(T,\mathbf{x})^T \\
	\end{array}
\right] \left[
	\begin{array}{ccc}
	    \mathbf{l}_{j+1}^+(T+\Delta T,\mathbf{x})\hspace{0.8mm}\cdots\hspace{0.8mm}\mathbf{l}_n^+(T+\Delta
T,\mathbf{x}) \\
	\end{array}
\right]\right\|_2 \\ &=\left\|\left[%
  \begin{aligned}
  &\langle \mathbf{l}_1^+(T,\mathbf{x}),\mathbf{l}_{j+1}^+(T+\Delta T,\mathbf{x}) \rangle
  \hspace{0.8mm}\cdots\hspace{0.8mm} \langle \mathbf{l}_1^+(T,\mathbf{x}),\mathbf{l}_n^+(T+\Delta
T,\mathbf{x})\rangle\\
  &\hspace{17mm}\vdots \hspace{40mm}\vdots \\	
  &\langle \mathbf{l}_j^+(T,\mathbf{x}),\mathbf{l}_{j+1}^+(T+\Delta T,\mathbf{x})
\rangle\hspace{0.8mm}\cdots\hspace{0.8mm} \langle
	    \mathbf{l}_j^+(T,\mathbf{x}),\mathbf{l}_n^+(T+\Delta T,\mathbf{x}) \rangle \\
  \end{aligned}%
\right]\right\|_2
\end{aligned}
\end{equation}

Using a result from \cite{Orszag}, we have for $T>0$ to
$1^{st}$-order in the time increment $\Delta T$
\begin{equation}\label{lode} \mathbf{l}^+_m(T+\Delta T)
=(1+c\Delta T)\mathbf{l}^+_m(T)+ \Delta T \sum_{i=1 (i\ne m)}^n
\frac{\left[(\mathbf{n}_i^+)^T
(A^T+A)\mathbf{n}_m^+\right]\mathbf{l}^+_i}
{e^{(\mu^+_m-\mu^+_i)T}-e^{(\mu^+_i-\mu^+_m)T}},
\end{equation} where $A=D\mathbf{f}(\mathbf{x})$ is the system matrix of the linearized
dynamics (\ref{lindyn}), $\mathbf{n}^+_i$ is a vector from the SVD
of the transition matrix $\Phi(T,\mathbf{x})$ as defined in Section
\ref{FTLE}, $c$ is a constant that is inconsequential in the following
developments and is thus left unspecified, the $\mathbf{x}$
dependence has been suppressed, and all exponents and vectors in the
summation on the right-hand-side are evaluated at $(T,\mathbf{x})$.
It follows that the inner products in (\ref{inners}) are
\begin{equation} \langle \mathbf{l}^+_k(T,\mathbf{x}),\mathbf{l}^+_m(T+\Delta T,\mathbf{x})\rangle=   \Delta
T\frac{\left[(\mathbf{n}_k^+)^T
(A^T+A)\mathbf{n}_m^+\right]}
{e^{(\mu^+_m-\mu^+_k)T}-e^{(\mu^+_k-\mu^+_m)T}}.
\end{equation} Because $k\in \{1,\dots,j\}$ and $m\in
\{j+1,\dots,n\}$, we have
$\exp[(\mu^+_k(T,\mathbf{x})-\mu^+_m(T,\mathbf{x}))T]\le
\exp[-{\Delta\mu}^+_j(\mathbf{x}) T]$. Let $\overline
a=\max_{x\in \mathcal{Y}}\max_{i\in
\{1,2,\dots,n\}}|\lambda_i(A^T+A)|$, the maximum eigenvalue
magnitude of $A^T+A$ over the set $\mathcal{Y}$. And let
$\alpha=\exp(-2{\Delta\mu}^+_j(\mathbf{x}) T_1)$ for some
$T_1 > t_s$. Then for $T\ge T_1 > 0$ we have
\begin{equation}|\langle \mathbf{l}^+_k(T,\mathbf{x}),\mathbf{l}^+_m(T+\Delta T,\mathbf{x})\rangle |\le
\frac{\overline a \Delta T}{1-\alpha}
e^{-{\Delta\mu}^+_j(\mathbf{x}) T}.
\end{equation} Upper-bounding
the 2-norm by the Frobenius norm and taking
$K=\sqrt{j(n-j)}$$\frac{\overline a \Delta T}{1-\alpha}$, the bound in the
theorem follows. This bound is conservative, due to the use of the
Frobenius norm, but it shows the exponential rate of convergence.
Using the bound (\ref{expbdd_A}), one can show that the sequence of
iterates is Cauchy. Moreover this is true for every sufficiently
small $\Delta T$. Because the space of $j$-dimensional subspaces in
$T_\mathbf{x}\mathbb{R}^n$, a Grassmannian, with the distance given in
(\ref{distance}) as the metric, is complete, we conclude that $
{\cal L}^+_j(T,\mathbf{x})$ approaches a fixed subspace. This
subspace is  ${\cal L}^+_j(\mathbf{\mathbf{x}})$ defined in Section
\ref{asympt}, because all vectors in it have exponents less than or equal to
$\mu^+_j(\mathbf{x})$ and one can show that any vector not in the
subspace must have a larger exponent. The proof for backward time is
similar. $\hfill\blacksquare$

\section{Averaging Time Determination}
\label{AppendixB} In
order to automate determining the averaging time for $\mathbf{x}_j$ for the calculation of the FTLVs, the averaging
time is iteratively increased, without
restricting the forward and backward averaging times to be the same. For the computataions in \ref{mckManifold}, for each pair $(x_1,\lambda_1)$, the value of
$(x_2,\lambda_2)$ approximating a point on the invariant center manifold is computed using an algorithm consisting of two nested
iteration loops with $i$
indicating the inner-loop iteration and $k$ the outer-loop iteration, with $i,k=0, 1, 2, \dots$. The variables and
iteration indices follow the format: $T^{(k)}_{fwd}$, $T^{(k)}_{bwd}$, $\mathbf{x}_j^{(i, k)}$, $x_2^{(i, k)}$, and
$\lambda_2^{(i, k)}$.

\begin{enumerate}
\item Initialization: Set $\mathbf{x}_j^{(0, 0)}=\mathbf{x}_j$ and $(x_2^{(0, 0)}, \lambda_2^{(0, 0)})$ to the values of
those coordinates in $\mathbf{x}_j^{(0, 0)}$. Set
$T^{(0)}_{fwd}=\overline{T}=0.5$ and $T^{(0)}_{bwd}=\overline{T}=0.5$.
\item Inner-loop iteration $i+1$ at outer iteration $k$: Calculate $\mathbf{l}^-_1(T^{(k)}_{bwd},\mathbf{x}_j^{(i,k)})$
and
$\mathbf{l}^+_4(T^{(k)}_{fwd},\mathbf{x}^{(i,k)})$
and determine the values of $x_2^{(i+1,k)}$ and
$\lambda_2^{(i+1,k)}$ that satisfy
 \begin{equation}\begin{array}{ccl}\label{eqn:MCK13_model}
\left\langle\mathbf{l}_1^-(T^{(k)}_{bwd},\mathbf{x}_j^{(i, k)}),\mathbf{f}(\mathbf{x}_j^{(i+1, k)})
\right\rangle = 0\\
\left\langle\mathbf{l}_4^+(T^{(k)}_{fwd},\mathbf{x}_j^{(i, k)}),\mathbf{f}(\mathbf{x}_j^{(i+1, k)}) \right\rangle = 0.
\end{array}\end{equation}
For this example the unknowns appear linearly; thus analytical solutions for $x_2^{(i+1,k)}$ and
$\lambda_2^{(i+1,k)}$ can be obtained.
Iterate until the inner-loop stopping criteria are met.
The stopping criteria consider the relative change in the dependent variables from the previous iteration
and $\theta^{(i+1,k)}$ is the angle between $\mathbf{f}(\mathbf{x}_j^{(i+1,k)})$ and its orthogonal projection in
$\mathcal{E}^c(T^{(k)}_{fwd},T^{(k)}_{bwd},\mathbf{x}_j^{(i,k)})$
 according to
\begin{equation}\begin{array}{ccl}
\label{eqn:MCK14_model}
\lvert x_2^{(i+1,k)}-x_2^{(i,k)}|/|x_2^{(i,k)}|<tol_{x_2} ,\\
|\lambda_2^{(i+1,k)}-\lambda_2^{(i,k)}|/|\lambda_2^{(i,k)}|<tol_{\lambda_2} ,\\
\theta^{(i+1,k)}<tol_{\theta} .\\
\end{array}
\end{equation}
 For this example, we used
$tol_{x_2}=tol_{\lambda_2}=tol_{\theta}=10^{-5}$. The approximation at the end of the inner-loop is denoted by
$\hat{\mathbf{x}}_j^{(k)}$.
\item Outer-loop iteration: Check the outer-loop stopping criterion
\begin{equation}\label{eqn:MCK26_model}
\|\hat{\mathbf{x}}_j^{(k)}-\hat{\mathbf{x}}_j^{(k-1)}\|_2<tol
\end{equation}
We used $tol=10^{-6}$. When $k=0$, we use $\mathbf{x}_j$ in place of $\hat{\mathbf{x}}_j^{(k-1)}$. If the criterion is
satisfied, stop and yield the final approximation $\hat{\mathbf{x}}_j$ to the center manifold point for the pair
$(x_1,\lambda_1)$ under consideration. Otherwise perform the $(k+1)^{th}$ outer-loop iteration with the averaging times
\begin{equation}
\label{eqn:MCK25_model}
T_{fwd}^{(k+1)} = T_{fwd}^{(k)}+dT_{fwd}, \hspace{5mm} T_{bwd}^{(k+1)} = T_{bwd}^{(k)}+dT_{bwd}. \\
\end{equation}
We used $dT_{fwd}=0.3$ and $dT_{bwd}=0.1$. With the new averaging times, repeat the inner-loop iterations starting with
$\hat{\mathbf{x}}_j^{(k)}$. \end{enumerate}
The computations for the five points required about $5$ inner iterations for each outer iteration and the forward and
backward averaging times were increased to about $5.0$ and $2.0$ respectively. Experiments with initializing the
iterative process with different dependent variable estimates consistently led to the same invariant center manifold point
approximations.

\end{document}